%% file: ms.tex
\documentclass[final, onefignum, onetabnum]{siamart171218}

\input{tex/ex_shared}

\ifpdf
\hypersetup{
  pdftitle={Dynamic Optimization of Thermodynamically Rigorous Models of Multiphase Flow in Porous Subsurface Oil Reservoirs},
  pdfauthor={Tobias K. S. Ritschel and John Bagterp J\O rgensen}
}
\fi

\begin{document}
	\maketitle
	
	\begin{abstract}
		\input{tex/abstract}
	\end{abstract}
	
	
	
	\input{tex/intro}
	\input{tex/thermal}
	\input{tex/isothermal}
	\input{tex/opt}
	\input{tex/ex}
	\input{tex/concl}
	
	\appendix
	
	\bibliographystyle{siamplain}
	

\input{build/ms.bbl}
	\input{tex/app_relperm}
	\input{tex/app_visc}
\end{document}

%% file: tex/ex_shared.tex
\usepackage{lipsum}
\usepackage{amsfonts}
\usepackage{graphicx}
\usepackage{epstopdf}
\usepackage{algorithmic}
\ifpdf
  \DeclareGraphicsExtensions{.eps,.pdf,.png,.jpg}
\else
  \DeclareGraphicsExtensions{.eps}
\fi

\newsiamremark{remark}{Remark}
\newsiamremark{hypothesis}{Hypothesis}
\crefname{hypothesis}{Hypothesis}{Hypotheses}
\newsiamthm{claim}{Claim}

\headers{Dynamic Optimization of Flow in Oil Reservoirs}{Tobias K. S. Ritschel and John Bagterp J\O rgensen}

\title{Dynamic Optimization of Thermodynamically Rigorous Models of Multiphase Flow in Porous Subsurface Oil Reservoirs%
	\thanks{
		\textbf{Funding:} This work was funded by Innovation Fund Denmark in the OPTION project (63-2013-3).
	} 
} 

\author{
	Tobias K. S. Ritschel\footnotemark[2]
	\and John Bagterp J\O rgensen\thanks{Department of Applied Mathematics and Computer Science \& Center for Energy Resources Engineering (CERE), Technical University of Denmark, DK-2800 Kgs. Lyngby, Denmark (\email{tobk@dtu.dk}, \email{jbjo@dtu.dk}).
	} 
} 

\usepackage{amsopn}


%% file: tex/abstract.tex
In this paper, we consider dynamic optimization of thermal and isothermal oil recovery processes which involve multicomponent three-phase flow in porous media. We present thermodynamically rigorous models of these processes based on 1) conservation of mass and energy, and 2) phase equilibrium. The conservation equations are partial differential equations. The phase equilibrium problems that are relevant to thermal and isothermal models are called the UV and the VT flash, and they are based on the second law of thermodynamics. We formulate these phase equilibrium problems as optimization problems and the phase equilibrium conditions as the corresponding first order optimality conditions. We demonstrate that the thermal and isothermal flow models are in a semi-explicit differential-algebraic form, and we solve the dynamic optimization problems with a previously developed gradient-based algorithm implemented in C/C++. We present numerical examples of optimized thermal and isothermal oil recovery strategies and discuss the computational performance of the dynamic optimization algorithm in these examples.

%% file: tex/keywords.tex
dynamic optimization,
single-shooting,
the adjoint method,
thermal and isothermal oil recovery,
multicomponent multiphase flow,
phase equilibrium,
UV flash,
VT flash

%% file: tex/classifications.tex
35Q35,
%
%
35Q93,
%
%
%
%
%
%
76N25,
%
%
76S05,
%
%
76T30,
%
%
93C20
%
%
%
%

%% file: tex/intro.tex
\section{Introduction}
Dynamic optimization is concerned with the computation of an optimal open-loop control strategy for a dynamical process. The objective of the optimization is either to 1) minimize the distance to predefined setpoints or 2) optimize the economics of the process. Dynamic optimization of multiphase flow processes in porous media is relevant to numerous engineering applications, e.g. production of oil from subsurface reservoirs \cite{Garipov:etal:2018, Kourounis:etal:2014, Polivka:Mikyska:2014, Ritschel:Jorgensen:2018c, Ritschel:Jorgensen:2018e, Zaydullin:etal:2014}, geothermal energy systems \cite{Thorvaldsson:Palsson:2012}, groundwater contamination and remediation \cite{Jang:Aral:2009}, trickle bed reactors \cite{Hannaoui:etal:2015}, fuel cells \cite{Berning:etal:2009, Kone:etal:2017}, food processing \cite{Khan:etal:2018}, and several others \cite{Pesavento:etal:2017}.

In this work, we consider dynamic optimization of thermal (varying temperature) and isothermal (constant temperature) oil recovery processes which involve multicomponent multiphase flow in porous rock. Oil recovery processes are described as primary, secondary, or tertiary \cite{Chen:2007, Chen:2006}. In primary recovery processes, the oil is recovered by means of the initial pressure in the reservoir. In secondary recovery processes, water is injected into the reservoir in order to maintain a high pressure. Tertiary recovery involves chemical, biological, or thermal injection with the purpose of mobilizing and recovering the oil that remains after the primary and secondary recovery processes. We are concerned with dynamic optimization of the secondary recovery process (also called waterflooding). However, dynamic optimization is equally applicable to the tertiary recovery processes (also called enhanced oil recovery processes). The objective of the dynamic optimization is to compute a field-wide production strategy that optimizes a long-term financial measure of the oil production, e.g. the total recovery or the net present value over the life-time of the oil reservoir.

Models of thermal and isothermal reservoir flow are based on two main principles: 1) conservation of mass and energy and 2) phase equilibrium. The conservation of energy is related to the first law of thermodynamics, while the equilibrium between phases is related to the second law of thermodynamics. The conservation equations are partial differential equations, and we formulate the phase equilibrium problems as inner optimization problems \cite{Michelsen:1999}. Consequently, the dynamic optimization problem that we consider belongs to the class of bilevel optimization problems \cite{Colson:etal:2007} as well as the closely related class of mathematical programs with equilibrium constraints \cite{Luo:etal:1996, Outrata:etal:1998}. We use the method of lines and discretize the conservation equations with a finite volume method in which the reservoir is represented by a discrete grid. The result of the discretization is a set of differential equations for each cell in the grid. Furthermore, we enforce the condition of phase equilibrium in each grid cell. The phase equilibrium conditions in the thermal model are different from the phase equilibrium conditions in the isothermal model. However, both sets of conditions are derived from the second law of thermodynamics which states that the entropy of a closed system in equilibrium is maximal \cite{Callen:1985, Koretsky:2013, Smith:etal:2005}. The phase equilibrium optimization problem in the thermal model is the UV (or UVn) flash which is a direct statement of the second law of thermodynamics. The objective of the UV flash optimization problem is to maximize entropy subject to constraints on the internal energy, $U$, the volume, $V$, and the total composition, $n$, i.e. total amount of moles of each component. The internal energy and the total composition are determined by the conservation equations while the volume is the size of the grid cell. The solution to the UV flash is the equilibrium temperature, pressure, and phase compositions. Isothermal systems are not closed. Consequently, the condition of maximal entropy does not apply directly. Instead, the Helmholtz energy is minimal for isothermal systems in equilibrium \cite{Callen:1985}. The phase equilibrium optimization problem in the isothermal model is the VT (or VTn) flash which involves minimization of the Helmholtz energy subject to constraints on the volume, $V$, the temperature, $T$, and the total composition, $n$. The solution to the VT flash is the equilibrium pressure and phase compositions.
In the reservoir simulation and optimization literature, it is common to formulate the phase equilibrium conditions as the isofugacity condition \cite{Chen:2007, Chen:2006, Garipov:etal:2018, Kourounis:etal:2014, Zaydullin:etal:2014} which is derived from the PT (or PTn) flash \cite{Ritschel:Jorgensen:2017, Ritschel:Jorgensen:2018b}. The PT flash is relevant to isothermal and isobaric (constant pressure) systems. For such systems, the Gibbs energy is minimal at equilibrium. Consequently, the PT flash involves minimization of the Gibbs energy subject to constraints on the temperature, $T$, pressure, $P$, and total composition, $n$. The solution to the PT flash is the equilibrium phase compositions. The UV and the VT flash are related to the PT flash \cite{Ritschel:etal:2017, Ritschel:Jorgensen:2018e}, and the condition of isofugacity in thermal and isothermal compositional reservoir flow models is derived from the UV and the VT flash.

Dynamic optimization of reservoir flow models most often involves models of immiscible two-phase flow \cite{Bukshtynov:etal:2015, Capolei:etal:2012, Codas:etal:2017, Heirung:etal:2011}, partially miscible two-phase flow \cite{Simon:Ulbrich:2015}, or polymer flooding \cite{Lei:etal:2012, Zhang:Li:2007}. Garipov et al. \cite{Garipov:etal:2018} and Zaydullin et al. \cite{Zaydullin:etal:2014} consider the simulation of thermal and compositional reservoir flow models, and Kourounis et al. \cite{Kourounis:etal:2014} present a gradient-based algorithm for dynamic optimization of isothermal and compositional reservoir flow models. However, none of the above models involve thermodynamically rigorous phase equilibrium conditions based on the UV or the VT flash. Pol\'{i}vka and Miky\v{s}ka \cite{Polivka:Mikyska:2014} consider simulation of an isothermal and compositional model that involves the VT flash.
Dynamic optimization of UV flash processes was first addressed by Ritschel et al. \cite{Ritschel:etal:2017, Ritschel:Capolei:Jorgensen:2017b, Ritschel:Capolei:Jorgensen:2017}. Furthermore, dynamic optimization of thermal and isothermal compositional reservoir flow models based on the UV and the VT flash was first considered by Ritschel and J{\o}rgensen \cite{Ritschel:Jorgensen:2018c, Ritschel:Jorgensen:2018e}.

There exists a number of algorithms for dynamic optimization of nonlinear systems \cite{Binder:etal:2001b}. Single-shooting algorithms involve the solution of numerical optimization problems in which the number of decision variables is independent of the number of state variables in the model. Therefore, single-shooting algorithms are often used for dynamic optimization of reservoir flow models \cite{Bukshtynov:etal:2015, Capolei:etal:2012, Forouzanfar:etal:2013} which typically involve a large number of state variables, i.e. on the order of $10^4$ - $10^7$ state variables. Alternative algorithms include multiple-shooting and simultaneous collocation which both require the solution of numerical optimization problems in which the number of decision variables does depend on the number of state variables. For reservoir models, the solution of such optimization problems can be intractable due to both high computation time and excessive memory requirements. However, both multiple-shooting \cite{Codas:etal:2017} and simultaneous collocation \cite{Heirung:etal:2011} have been used for dynamic optimization of reservoir flow models. Efficient algorithms for the solution of numerical optimization problems require the gradients of the objective function. For single- and multiple-shooting algorithms, such gradients can be computed efficiently with either an adjoint method \cite{Capolei:etal:2012, Jorgensen:2007b, Volcker:etal:2011} or a forward method \cite{Kristensen:etal:2005}. Alternatives to gradient-based optimization algorithms include stochastic approximation methods \cite{Zhao:etal:2016} and metaheuristic methods \cite{Onwunalu:Durlofsky:2010}.

In this work, we present thermodynamically rigorous models of thermal and isothermal waterflooding processes. We use the method of lines and discretize the conservation equations with a finite volume method. We demonstrate that the resulting equations are in a semi-explicit index-1 differential-algebraic form. Ritschel et al. \cite{Ritschel:etal:2017} describe a gradient-based dynamic optimization algorithm for such systems. The algorithm uses a single-shooting method together with an adjoint method for the computation of gradients. We implement the algorithm in C/C++ based on the open-source software DUNE \cite{Bastian:etal:2008, Bastian:etal:2008b, Blatt:Bastian:2007}, the open-source software ThermoLib \cite{Ritschel:etal:2016, Ritschel:etal:2017b}, and the commercial software KNITRO (IPOPT \cite{Wachter:Biegler:2006} is an open-source alternative to KNITRO). We use the thermodynamic software ThermoLib to evaluate thermodynamic properties based on the Peng-Robinson equation of state. The ThermoLib routines also evaluate the first and second order derivatives of the thermodynamic functions which are necessary in the gradient-based dynamic optimization algorithm. Finally, we present numerical examples of optimized thermal and isothermal waterflooding strategies, and we discuss the computational performance of the C/C++ implementation.

The remainder of this paper is organized as follows. We present the thermal and compositional model in Section \ref{sec:thermal}, and we present the isothermal and compositional model in Section \ref{sec:isothermal}. We formulate the dynamic optimization problem and discuss the C/C++ implementation in Section \ref{sec:opt}. In Section \ref{sec:ex}, we present the numerical examples, and we present conclusions in Section \ref{sec:concl}.

%% file: tex/thermal.tex
\section{Thermal and compositional reservoir flow model}\label{sec:thermal}
\begin{figure}[t]
	\centering
	\includegraphics[width=0.8\textwidth]{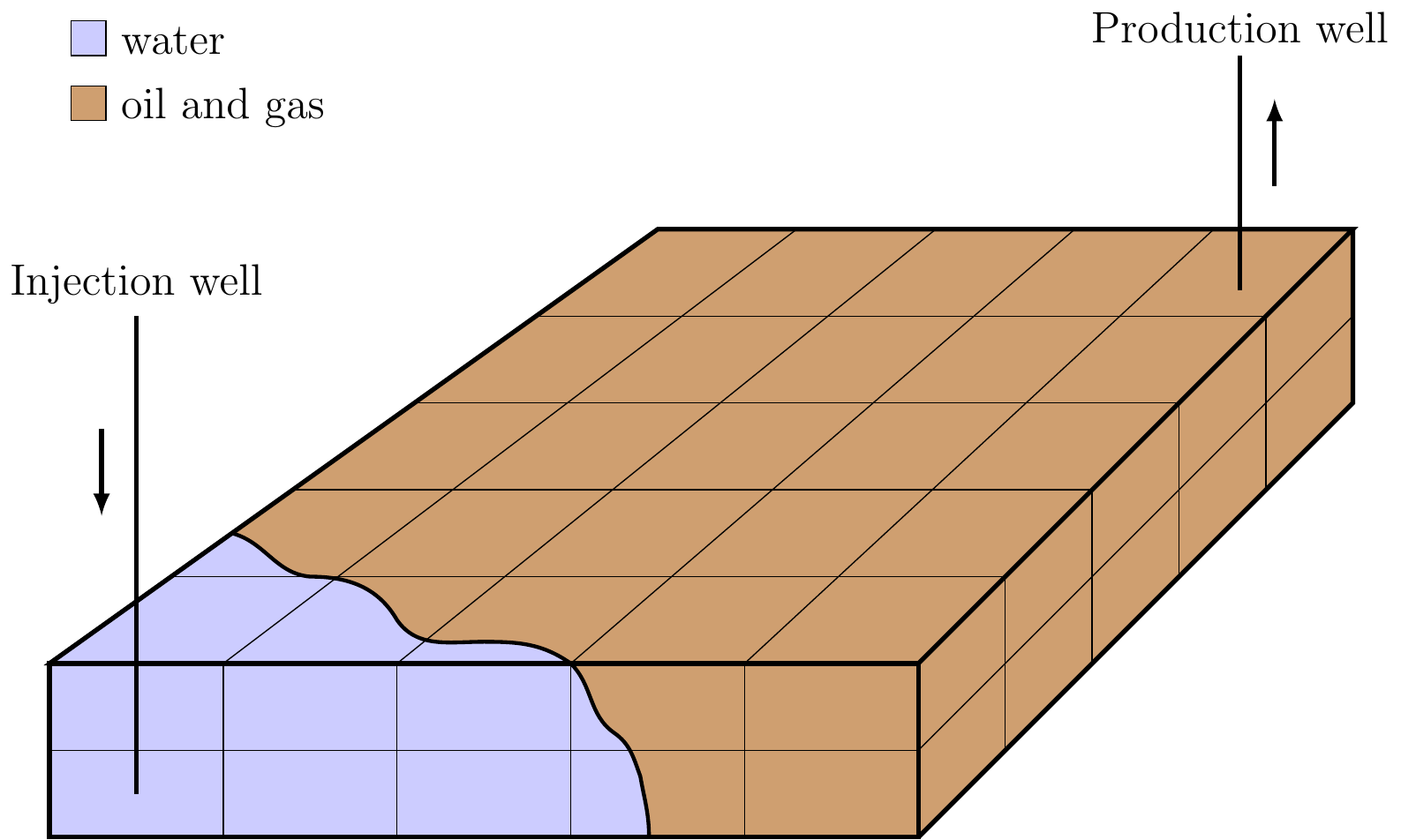}
	\caption{Illustration of the waterflooding process.}
	\label{fig:res}
\end{figure}

In this section, we describe the thermal and compositional model. The waterflooding process is illustrated in Fig. \ref{fig:res} for a rectangular (and discretized) reservoir. The model consists of a set of mass conservation equations, one energy balance equation, and a set of phase equilibrium conditions. The flow of mass in the reservoir is due to advection while the flow of energy is due to both advection and conduction. The phase equilibrium problem is the UV flash, and we assume that the fluid and the rock reach thermal, mechanical, and chemical equilibrium instantaneously, i.e. that they are in equilibrium at all times. We use a finite volume method to discretize the conservation equations, and we enforce the phase equilibrium in each cell of the discretized reservoir. Finally, we demonstrate that the thermal and compositional model is in the semi-explicit index-1 differential-algebraic form.

%
%

\subsection{Phase equilibrium}\label{sec:uv}
\input{tex/uv}

\subsection{Conservation of mass}\label{sec:mass}
\input{tex/mass}

\subsection{Conservation of energy}\label{sec:energy}
\input{tex/energy}

\subsection{Well equations}\label{sec:wells}
\input{tex/wells}

\subsection{Darcy's law}\label{sec:darcy}
\input{tex/darcy}

\subsection{Relative permeability}\label{sec:relperm}
\input{tex/relperm}

\subsection{Viscosity}\label{sec:visc}
\input{tex/visc}

\subsection{Thermodynamics}\label{sec:thermo}
\input{tex/thermo}

\subsection{Finite volume discretization}\label{sec:fv}
\input{tex/disc}

\subsection{The semi-explicit differential-algebraic form}\label{sec:sedae}
\input{tex/sedae}

%% file: tex/uv.tex
\begin{figure}[t]
	\centering
	\includegraphics[width=0.85\textwidth]{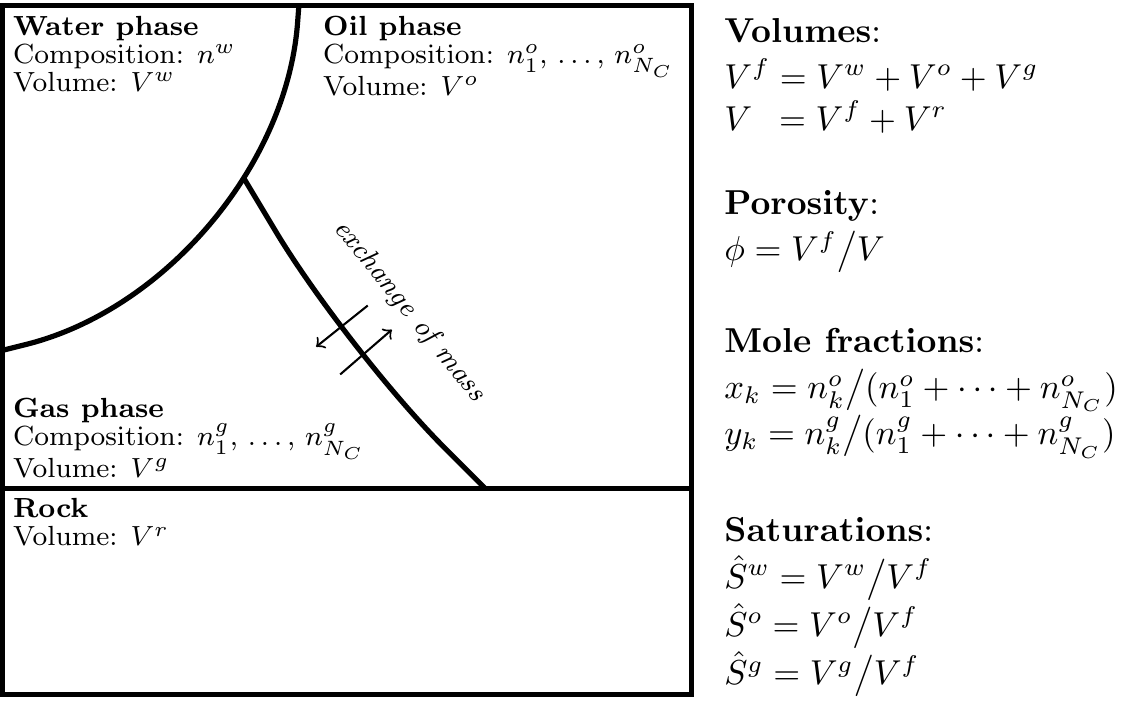}
	\caption{Illustration of the fluid phases and the rock in each grid cell.}
	\label{fig:cv}
\end{figure}

Each grid cell in the discretized reservoir contains a water phase ($w$), an oil phase ($o$), a gas phase ($g$), and a (solid) rock phase ($r$) as illustrated in Fig. \ref{fig:cv}.
The water phase is immiscible with the oil and the gas phase. The water phase contains only water while the oil and the gas phase contain $N_C$ chemical components. The fluid phases and the rock are in thermal and mechanical equilibrium, i.e. $T^\alpha = T$ and $P^\alpha = P$ for $\alpha\in\{w, o, g, r\}$. Furthermore, the oil and the gas phase are in chemical equilibrium. The UV flash optimization problem describing the isoenergetic-isochoric chemical equilibrium is
\begin{subequations}\label{eq:model:uv}
	\begin{align}
		\max_{T, P, n^w, n^o, n^g}	\qquad	& S^w + S^o + S^g + S^r, \\
		\textrm{subject to}			\qquad	& U^w + U^o + U^g + U^r = U, \\
		& V^w + V^o + V^g + V^r = V, \\
		& n^w = n_w, \\
		& n_k^o + n_k^g = n_k, \quad k = 1, \ldots, N_C.
	\end{align}
\end{subequations}
$S^\alpha = S^\alpha(T, P, n^\alpha)$, $U^\alpha = U^\alpha(T, P, n^\alpha)$, and $V^\alpha = V^\alpha(T, P, n^\alpha)$ are the entropy, internal energy, and volume of phase $\alpha\in\{w, o, g, r\}$, respectively. $U$ and $V$ are the specified internal energy and volume. $n_w$ and $n_k$ are the specified total amount of moles of water and of component $k$. $n_w$, $n_k$, and $U$ are determined by the mass and energy conservation equations, and $V$ is the size of the grid cell in the discretized reservoir. The solution to \eqref{eq:model:uv} is the equilibrium temperature, $T$, pressure, $P$, and phase compositions, $n^\alpha$ for $\alpha\in\{w, o, g\}$.

%% file: tex/mass.tex
Each mass conservation equation contains 1) a flux term related to the flow in the reservoir and 2) a source term related to the injection of water and the production of the reservoir fluid:
\begin{subequations}\label{eq:model:mass}
	\begin{align}
		\partial_t C_w &= -\nabla\cdot\textbf{N}^w + Q^w, \\
		\partial_t C_k &= -\nabla\cdot\textbf{N}_k + Q_k, & k &= 1, \ldots, N_C.
	\end{align}
\end{subequations}
$C_w$ and $C_k$ are the molar concentrations of water and component $k$. $\mathbf{N}^w$ is the molar flux of the water phase, and the molar flux of the $k$'th component is
\begin{align}
	\textbf{N}_k &= x_k \textbf{N}^o + y_k \textbf{N}^g, & k &= 1, \ldots, N_C.
\end{align}
$\mathbf{N}^\alpha$ is the molar flux of phase $\alpha\in\{o, g\}$. $x_k$ and $y_k$ are the oil and gas mole fractions of component $k$, i.e the moles of component $k$ in the oil and the gas phase divided by the total amount of moles in the respective phase. The source terms describe the molar well flow rates:
\begin{subequations}
	\begin{align}
		Q^w &= Q^{w, \text{inj}} - Q^{w, \text{prod}}, \\
		Q_k &= -\left(x_k Q^{o, \text{prod}} + y_k Q^{g, \text{prod}}\right), & k &= 1, \ldots, N_C.
	\end{align}
\end{subequations}
In Section \ref{sec:wells}, we provide expressions for the molar injection rate of the water phase, $Q^{w, \text{inj}}$, and the molar production rate of phase $\alpha\in\{w, o, g\}$, $Q^{\alpha, \text{prod}}$.

%% file: tex/energy.tex
First, we describe the conservation of energy of the fluid ($f$) and the rock ($r$) separately, i.e. without assuming thermal equilibrium between the fluid and the rock. Consequently, we distinguish between the temperature of the fluid phases, $T^f$, and the temperature of the rock, $T^r$. Next, we incorporate the assumption of thermal equilibrium and present the energy conservation equation for the combined rock-fluid system. The energy conservation equations for the fluid and the rock are
\begin{subequations}
	\begin{align}
		\partial_t u^f &= -\nabla\cdot\textbf{N}_u^f + Q_u^f,	\label{eq:model:energy:fluid} \\
		\partial_t u^r &= -\nabla\cdot\textbf{N}_u^r + Q_u^r. 	\label{eq:model:energy:rock}
	\end{align}
\end{subequations}
$u^f$ and $u^r$ are the internal energies per unit volume of the fluid and the rock. The heat flux of the fluid is due to the advective flow of the three fluid phases:
\begin{align}\label{eq:model:energy:fluid:flux}
	\textbf{N}_u^f &= h^w \textbf{N}^w + h^o \textbf{N}^o + h^g \textbf{N}^g.
\end{align}
$h^\alpha = h^\alpha(T, P, n^\alpha)$ is the molar enthalpy of phase $\alpha\in\{w, o, g\}$. The heat flux of the rock is due to conduction, and we describe it using Fourier's law of thermal conduction \cite[Chap.~1]{Holman:2010}:
\begin{align}\label{eq:model:energy:rock:flux}
	\textbf{N}_u^r &= -k_T^r\nabla T^r.
\end{align}
$k_T^r$ is the thermal conductivity of the rock. The source term in the fluid energy balance \eqref{eq:model:energy:fluid} describes 1) the transfer of energy through the wells and 2) the transfer of energy through the rock-fluid interface:
\begin{align}\label{eq:model:energy:fluid:source}
	Q_u^f &= h^{w, \text{inj}} Q^{w, \text{inj}} - \sum_{\alpha\in\{w, o, g\}} h^\alpha Q^{\alpha, \text{prod}} + Q^{rf}.
\end{align}
$h^{w, \text{inj}}$ is the molar enthalpy of the injected water. We describe the thermal conduction through the rock-fluid interface using Newton's law of cooling \cite[Chap.~1]{Holman:2010}:
\begin{align}
	Q^{rf} &= -k_T^{rf}(T^f - T^r).
\end{align}
$k_T^{rf}$ is the thermal conductivity of the rock-fluid interface. The source term in the energy balance for the rock contains terms describing 1) the transfer of energy through the rock-fluid interface and 2) the transfer of energy to the surroundings ($s$) of the reservoir:
\begin{align}
	Q_u^r &= -Q^{rf} - Q^{rs}.
\end{align}
We describe the thermal conduction through the interface between the rock and the surroundings of the reservoir using Newton's law of cooling:
\begin{align}
	Q^{rs} &= -k_T^{rs}(T^s - T^r).
\end{align}
$T^s$ is the temperature of the surroundings, and $k_T^{rs}$ is the thermal conductivity of the interface between the rock and surroundings. Now, we assume that energy is transferred instantaneously between the fluid and the rock, i.e. that the thermal conductivity of the rock-fluid interface, $k_T^{rf}$, is infinite. Consequently, the temperature of the fluid and the rock are equal, i.e. $T^f = T^r = T$. In order to obtain a conservation equation for the internal energy of the combined rock-fluid system, $u = u^f + u^r$, we add \eqref{eq:model:energy:fluid} and \eqref{eq:model:energy:rock}:
\begin{align}\label{eq:model:energy:tot}
	\partial_t u &= -\nabla\cdot\textbf{N}_u + Q_u.
\end{align}
The heat flux, $\mathbf{N}_u$, and the source term, $Q_u$, are
\begin{subequations}
	\begin{align}
		\mathbf{N}_u &= h^w \textbf{N}^w + h^o \textbf{N}^o + h^g \textbf{N}^g - k_T^r\nabla T, \\
		Q_u &= h^{w, \text{inj}}Q^{w, \text{inj}} - \sum_{\alpha\in\{w, o, g\}} h^\alpha Q^{\alpha, \text{prod}} - Q^{rs}.
	\end{align}
\end{subequations}

%% file: tex/wells.tex
The injection and the production wells are perforated in certain places in the reservoir, i.e. the injection and production source terms will only be nonzero in a few locations. The model of the well flow depends on the discretization of the reservoir. For a given grid cell, the molar injection and production phase flow rates are
\begin{subequations}
	\begin{align}
	Q^{w, 		\text{inj }} &= \frac{1}{V} \text{WI} \rho^w      \frac{k_r^w     }{\mu^w     }\left(P^\text{bhp} - P      \right), \\
	Q^{\alpha, 	\text{prod}} &= \frac{1}{V} \text{WI} \rho^\alpha \frac{k_r^\alpha}{\mu^\alpha}\left(P       - P^\text{bhp}\right), & \alpha &\in \{w, o, g\}.
	\end{align}
\end{subequations}
$V$ is the volume of the perforated grid cell, and WI is the well index which is a scalar quantity that describes the ability of the perforation to transmit fluid. $\rho^\alpha = \rho^\alpha(T, P, n^\alpha)$, $k_r^\alpha = k_r^\alpha(T, P, n^w, n^o, n^g)$, and $\mu^\alpha = \mu^\alpha(T, P, n^\alpha)$ are the molar density, the relative permeability, and the viscosity of phase $\alpha\in\{w, o, g\}$. $P^\text{bhp}$ is the bottom-hole pressure (BHP) in the well.

%% file: tex/darcy.tex
The molar phase flux is the product of the molar density and the volumetric phase flux:
\begin{align}
	\mathbf{N}^\alpha &= \rho^\alpha \mathbf{u}^\alpha, & \alpha &\in \{w, o, g\}.
\end{align}
We describe the volumetric phase flux with Darcy's law:
\begin{align}
\mathbf{u}^\alpha &= - \frac{k_r^\alpha}{\mu^\alpha} \mathbf{K} \left(\nabla P - \rho^\alpha g \nabla z\right), & \alpha &\in \{w, o, g\}.
\end{align}
\textbf{K} is a permeability tensor, $g$ is the gravity acceleration, and $z$ is the depth.

%% file: tex/relperm.tex
We use Stone's model II to describe the relative permeabilities \cite{Delshad:Pope:1989}. The relative permeabilities depend on the phase saturations, $\hat{S}^\alpha = V^\alpha/(V^w + V^o + V^g)$ for $\alpha\in\{w, o, g\}$. Consequently, the relative permeability of phase $\alpha\in\{w, o, g\}$ depends on the temperature, pressure, and the compositions of all phases:
\begin{align}
k_r^\alpha &= k_r^\alpha(T, P, n^w, n^o, n^g), & \alpha &\in \{w, o, g\}.
\end{align}
In Appendix \ref{app:relperm}, we present the expressions for the relative permeabilities in detail.

%% file: tex/visc.tex
We describe the viscosities of the oil and the gas phase with the model by Lohrenz et al. \cite{Lohrenz:etal:1964}, and we model the water viscosity by $(1/\mu^w)(\partial \mu^w/\partial P) = c_\mu^w$ where $c_\mu^w$ is the viscosibility of the water phase. Consequently, the viscosities are functions of temperature, pressure, and the phase compositions:
\begin{align}
\mu^\alpha &= \mu^\alpha(T, P, n^\alpha), & \alpha &\in \{w, o, g\}.
\end{align}
In Appendix \ref{app:visc}, we describe the viscosity of the oil and gas phases in detail.

%% file: tex/thermo.tex
The phase equilibrium optimization problem \eqref{eq:model:uv}, the fluid heat flux \eqref{eq:model:energy:fluid:flux}, and the fluid heat source \eqref{eq:model:energy:fluid:source} involve thermodynamical functions. We use the open-source thermodynamic software ThermoLib \cite{Ritschel:etal:2016, Ritschel:etal:2017b} to evaluate the enthalpy, entropy, and volume of the fluid phases:
%
%
\begin{subequations}
	\begin{align}
	H^\alpha &= H^\alpha(T, P, n^\alpha), & \alpha &\in \{w, o, g\}, \\
	S^\alpha &= S^\alpha(T, P, n^\alpha), & \alpha &\in \{w, o, g\}, \\
	V^\alpha &= V^\alpha(T, P, n^\alpha), & \alpha &\in \{w, o, g\}.
	\end{align}
\end{subequations}
The thermodynamic model in ThermoLib is based on data and correlations from the DIPPR database \cite{Thomson:1996} as well as cubic equations of state \cite{Gmehling:etal:2012, Koretsky:2013, Smith:etal:2005}. We use the Peng-Robinson equation of state \cite{Peng:robinson:1976}.
The first order optimality conditions of the phase equilibrium problem \eqref{eq:model:uv} are algebraic equations. The gradient-based dynamic optimization algorithm described by Ritschel et al. \cite{Ritschel:etal:2017} requires the Jacobian matrices of these algebraic equations. Consequently, the algorithm requires both the first and second order derivatives of the thermodynamic functions with respect to temperature, pressure, and mole numbers. The ThermoLib routines evaluate such derivatives based on the analytical expressions described by Ritschel et al. \cite{Ritschel:etal:2016}.
The thermodynamic properties of the rock, $H^r = H^r(T, P)$, $S^r = S^r(T, P)$, and $V^r = V^r(T, P)$ are also computed from an equation of state. We use a temperature-independent equation of state, $(1/V^r)(\partial V^r/\partial P) = c^r$, and we assume that the rock compressibility, $c^r$, is constant.
We compute other thermodynamic functions with the fundamental thermodynamic relations $U^\alpha = H^\alpha - PV^\alpha$, $G^\alpha = H^\alpha - TS^\alpha$, and $A^\alpha = U^\alpha - TS^\alpha$ for $\alpha\in\{w, o, g, r\}$.

%% file: tex/disc.tex
The mass conservation equations \eqref{eq:model:mass} and the energy conservation equation \eqref{eq:model:energy:tot} are in the form
\begin{align}\label{eq:disc:cons}
	\partial_t C = -\nabla\cdot\textbf{N} + Q.
\end{align}
In this section, we describe the finite volume discretization of \eqref{eq:disc:cons} and use it to discretize the mass and energy conservation equations. We consider a discretized reservoir that consists of a set of grid cells, $\{\Omega_i\}_{i\in\mathcal{N}}$, where $\mathcal{N}$ is a set of grid cell indices. We assume that each grid cell is a polyhedron and that each face of the polyhedron is shared by exactly two cells. We integrate \eqref{eq:disc:cons} over each of the grid cells and interchange integration and differentiation on the left-hand side:
\begin{align}\label{eq:fv}
	\partial_t \int_{\Omega_i} C \,dV &= -\int_{\Omega_i} \nabla\cdot\textbf{N} \,dV + \int_{\Omega_i} Q \,dV, & i &\in \mathcal{N}.
\end{align}
We apply Gauss' divergence theorem to the first integral on the right-hand side:
\begin{align}\label{eq:disc:gauss}
	\int_{\Omega_i} \nabla\cdot\textbf{N} \,dV &= \int_{\partial\Omega_i} \textbf{N}\cdot\textbf{n} \,dA, & i &\in \mathcal{N}.
\end{align}
$\partial\Omega_i$ is the boundary of the $i$'th grid cell, and $\mathbf{n}$ is the outward normal vector. We split up the boundary integral over each of the faces of the cell:
\begin{align}\label{eq:disc:face}
	\int_{\partial\Omega_i} \textbf{N}\cdot\textbf{n} \,dA &= \sum_{j\in\mathcal{N}^{(i)}} \int_{\gamma_{ij}} \textbf{N}\cdot\textbf{n} \,dA, & i &\in \mathcal{N}.
\end{align}
$\mathcal{N}^{(i)}$ is the set of cells that share a face with the $i$'th grid cell, and $\gamma_{ij}$ is the face that is shared by the $i$'th and the $j$'th grid cell. We use quadrature to approximate the integral of the source term in \eqref{eq:fv} and the integral over $\gamma_{ij}$ in \eqref{eq:disc:face}:
\begin{subequations}
	\begin{align}
	\int_{\Omega_i} Q \,dV &\approx (QV)_i, & i &\in \mathcal{N}, \\
	\int_{\gamma_{ij}} \textbf{N}\cdot\textbf{n} \,dA &\approx (A\textbf{N}\cdot\textbf{n})_{ij}, & i &\in \mathcal{N}, & j &\in \mathcal{N}^{(i)}. \label{eq:model:fv:flux}
	\end{align}
\end{subequations}
The subscript $i$ indicates that a quantity is related to the $i$'th grid cell while the subscript $ij$ indicates that it is related to the face $\gamma_{ij}$. $V_i$ is the volume of $\Omega_i$, and $A_{ij}$ is the area of $\gamma_{ij}$. We now apply the finite volume discretization to the mass and energy conservation equations. The integrals of the internal energy per unit volume and the concentrations are evaluated exactly:
\begin{subequations}
	\begin{align}
		\int_{\Omega_i}   u \,dV &= U_{   i}, & i &\in \mathcal{N}, \\
		\int_{\Omega_i} C_w \,dV &= n_{w, i}, & i &\in \mathcal{N}, \\
		\int_{\Omega_i} C_k \,dV &= n_{k, i}, & i &\in \mathcal{N}.
	\end{align}
\end{subequations}
The right-hand side of \eqref{eq:model:fv:flux} involves the flux evaluated at the center of the face which we approximate with a two-point flux approximation \cite{Lie:2014}. The resulting approximation of the right-hand side of \eqref{eq:model:fv:flux} for the heat and mass fluxes are
\begin{subequations}\label{eq:disc:flux}
\begin{align}
	(A\mathbf{N}_u\cdot\mathbf{n})_{ij} &\approx -\sum_{\alpha\in\{w, o, g\}} (h^\alpha \Gamma \hat{H}^\alpha \Delta\Phi^\alpha)_{ij} + (\Gamma_T \Delta T)_{ij}, & i &\in \mathcal{N}, & j &\in \mathcal{N}^{(i)},
	\label{eq:disc:flux:u} \\
	(A\mathbf{N}^w\cdot\mathbf{n})_{ij} &\approx -(\Gamma \hat{H}^w \Delta\Phi^w)_{ij}, & i &\in \mathcal{N}, & j &\in \mathcal{N}^{(i)},
	\label{eq:disc:flux:w} \\
	(A\mathbf{N}_k\cdot\mathbf{n})_{ij} &\approx -(x_k \Gamma \hat{H}^o \Delta\Phi^o + y_k \Gamma \hat{H}^g \Delta\Phi^g)_{ij}, & i &\in \mathcal{N}, & j &\in \mathcal{N}^{(i)}.
	\label{eq:disc:flux:k}
\end{align}
\end{subequations}
The difference in temperature is $\Delta T_{ij} = T_j - T_i$. $\Gamma_{ij}$ is the geometric part of the transmissibilities:
\begin{subequations}\label{eq:disc:geom}
\begin{align}
	\Gamma_{ij} &= A_{ij}\left(\hat{\Gamma}_{ij}^{-1} + \hat{\Gamma}_{ji}^{-1}\right)^{-1}, & i &\in \mathcal{N}, & j &\in \mathcal{N}^{(i)}, \\
	\hat{\Gamma}_{ij} &= \left(\textbf{K}_i\frac{c_{ij} - c_i}{|c_{ij} - c_i|^2}\right)\cdot\textbf{n}_{ij}, & i &\in \mathcal{N}, & j &\in \mathcal{N}^{(i)}.
\end{align}
\end{subequations}
$c_{ij}$ is the center of $\gamma_{ij}$, $c_i$ is the center of $\Omega_i$, and $\hat{\Gamma}_{ij}$ is the one-sided transmissibility. The expression for $\Gamma_{T, ij}$ is analogous to \eqref{eq:disc:geom}. However, the thermal conductivity of the rock replaces the permeability tensor:
\begin{subequations}\label{eq:disc:geom:heat}
	\begin{align}
		\Gamma_{T, ij} &= A_{ij}\left(\hat{\Gamma}_{T, ij}^{-1} + \hat{\Gamma}_{T, ji}^{-1}\right)^{-1}, & i &\in \mathcal{N}, & j &\in \mathcal{N}^{(i)}, \\
		\hat{\Gamma}_{T, ij} &= \left(k_{T, i}^r\frac{c_{ij} - c_i}{|c_{ij} - c_i|^2}\right)\cdot\textbf{n}_{ij}, & i &\in \mathcal{N}, & j &\in \mathcal{N}^{(i)}.
	\end{align}
\end{subequations}
The difference in potential and the fluid part of the transmissibilities are
\begin{subequations}
\begin{align}
	\Delta \Phi_{ij}^\alpha &= \left(\Delta P - \rho^\alpha g \Delta z\right)_{ij}, & i &\in \mathcal{N}, & j &\in \mathcal{N}^{(i)}, \label{eq:model:fv:pot} \\
	\hat{H}_{ij}^\alpha &= 
	\begin{cases}
		(\rho^\alpha k_r^\alpha/\mu^\alpha)_i, & \Delta \Phi_{ij}^\alpha < 0, \\
		(\rho^\alpha k_r^\alpha/\mu^\alpha)_j, & \Delta \Phi_{ij}^\alpha \geq 0,
	\end{cases} & i &\in \mathcal{N}, & j &\in \mathcal{N}^{(i)}. \label{eq:disc:fv:fluidtrans}
\end{align}
\end{subequations}
The differences in pressure and depth are $\Delta P_{ij} = P_j - P_i$ and $\Delta z_{ij} = z_j - z_i$. We approximate the density on the face center by $\rho_{ij}^\alpha \approx (\rho_i^\alpha + \rho_j^\alpha)/2$. In \eqref{eq:disc:fv:fluidtrans}, we have upwinded the fluid part of the transmissibilities to ensure numerical stability. Similarly, we upwind $h^\alpha$ in \eqref{eq:disc:flux:u} as well as $x_k$ and $y_k$ in \eqref{eq:disc:flux:k}. The differential equations that result from the finite volume discretization of the mass and energy conservation equations are
\begin{subequations}\label{eq:cons:disc}
	\begin{align}
	\dot{U}_{   i} &= \sum_{j\in\mathcal{N}^{(i)}} \left(\sum_{\alpha\in\{w, o, g\}} (h^\alpha \Gamma \hat{H}^\alpha \Delta\Phi^\alpha)_{ij} + (\Gamma_T \Delta T)_{ij}\right) + (Q_u V)_i, & i &\in \mathcal{N}, \\
	\dot{n}_{w, i} &= \sum_{j\in\mathcal{N}^{(i)}} (\Gamma \hat{H}^w \Delta\Phi^w)_{ij} + (Q^w V)_i, & i &\in \mathcal{N}, \\
	\dot{n}_{k, i} &= \sum_{j\in\mathcal{N}^{(i)}} (x_k \Gamma \hat{H}^o \Delta\Phi^o + y_k \Gamma \hat{H}^g \Delta\Phi^g)_{ij} + (Q_k V)_i, & i &\in \mathcal{N}.
	\end{align}
\end{subequations}
The internal energy and the total amounts of moles on the left-hand side of \eqref{eq:cons:disc} appear as specified quantities in the phase equilibrium problem described in Section \ref{sec:uv}.

%% file: tex/sedae.tex
The phase equilibrium problem described in Section \ref{sec:uv} is in the form
\begin{subequations}\label{eq:opt:flash}
	\begin{align}
		\min_{y_i}  		\qquad 	& f(y_i), \label{eq:opt:flash:f} \\
		\text{subject to} 	\quad 	& g(y_i) = x_i, \label{eq:opt:flash:g} \\
		& h(y_i) = 0. \label{eq:opt:flash:h}
	\end{align}
\end{subequations}
$x_i = \begin{bmatrix} U; n_w; n \end{bmatrix}_i \in \mathbb{R}^{2 + N_C}$ is the state vector, and $y_i = \begin{bmatrix} T; P; n^w; n^o; n^g \end{bmatrix}_i \in \mathbb{R}^{3 + 2 N_C}$ is a vector of algebraic variables. The phase equilibrium conditions are the first order optimality conditions (also called Karush-Kuhn-Tucker or KKT conditions) of \eqref{eq:opt:flash}. The first order optimality conditions are a set of algebraic equations, $G_i(x_i, y_i, z_i) = 0$, because \eqref{eq:opt:flash} does not contain inequality constraints \cite{Nocedal:Wright:2006}. $z_i \in \mathbb{R}^{3 + N_C}$ are Lagrange multipliers.

The left-hand sides of the differential equations \eqref{eq:cons:disc} are the time derivatives of the state variables while the quantities on the right-hand sides depend exclusively on the algebraic variables in the $i$'th cell, $y_i$, and in neighbouring cells, $\{y_j\}_{j\in\mathcal{N}^{(i)}}$, as well as the manipulated inputs, $u_i = P_i^\text{bhp} \in \mathbb{R}$, and the disturbance variables, $d_i = T_i^\text{inj}\in\mathbb{R}$. The temperature of the injected water, $T^\text{inj}$, is used to evaluate the enthalpy of the injected water, $h^{w, \text{inj}}$. The manipulated inputs and the disturbance variables are only nonempty for cells that are perforated by a well. The differential equations \eqref{eq:cons:disc} are thus in the form $\dot{x}_i(t) = F_i(y_i(t), \{y_j(t)\}_{j\in\mathcal{N}^{(i)}}, u_i(t), d_i(t))$. Consequently, the collection of the differential equations and phase equilibrium conditions for all grid cells is in the semi-explicit differential-algebraic form,
\begin{subequations}\label{eq:sedae}
	\begin{align}
		&G(x(t), y(t), z(t)) = 0, \\
		&\dot{x}(t) = F(y(t), u(t), d(t)).
	\end{align}
\end{subequations}
Furthermore, the algebraic equations are of index 1, i.e. $G_i(x_i, y_i, z_i) = 0$ can be solved for $y_i$ and $z_i$ when $x_i$ is specified.

%% file: tex/isothermal.tex
\section{Isothermal and compositional reservoir flow model}\label{sec:isothermal}
In this section, we adapt the thermal and compositional model presented in Section \ref{sec:thermal} to isothermal systems. In isothermal systems, all involved thermal conductivities are infinite such that energy is transferred instantaneously between 1) the fluid and the rock, and 2) the rock and the surroundings, until thermal equilibrium is reached.
Furthermore, the heat capacity of the surroundings is infinite such that their temperature is constant despite the supply or removal of energy.
The key difference between the thermal and the isothermal model is the phase equilibrium problem which, for isothermal systems, is the VT flash. Furthermore, the isothermal model does not involve an energy conservation equation. However, the mass conservation equations in the two models are identical. Therefore, we only discuss 1) the phase equilibrium problem and 2) the semi-explicit differential-algebraic form of the model.

\subsection{Phase equilibrium}
\input{tex/vt}

\subsection{The semi-explicit differential-algebraic form}
\input{tex/isosedae}

%% file: tex/vt.tex
The VT flash optimization problem describing the isochoric-isothermal chemical equilibrium is
\begin{subequations}\label{eq:vt}
	\begin{align}
		\min_{P, n^w, n^o, n^g} \quad 	& A^w + A^o + A^g + A^r, \label{eq:vt:a} \\
		\text{subject to} 		\quad 	& V^w + V^o + V^g + V^r = V, \label{eq:vt:v} \\
										& n^w = n_w, \label{eq:vt:w} \\
										& n_k^o + n_k^g = n_k, & k &= 1, \ldots, N_C. \label{eq:vt:k}
	\end{align}
\end{subequations}
$A^\alpha = A^\alpha(T, P, n^\alpha)$  is the Helmholtz energy of phase $\alpha\in\{w, o, g, r\}$. The main differences between the VT flash \eqref{eq:vt} and the UV flash \eqref{eq:model:uv} are that in the VT flash 1) the Helmholtz energy is minimized, 2) there is no constraint on the internal energy, and 3) temperature is not an optimization variable.
%

%% file: tex/isosedae.tex
The VT flash optimization problem is in the same form as the UV flash optimization problem, i.e. \eqref{eq:opt:flash}. The state variables are $x_i = \begin{bmatrix} n_w; n \end{bmatrix}_i \in \mathbb{R}^{1 + N_C}$, and the algebraic variables are $y_i = \begin{bmatrix} P; n^w; n^o; n^g \end{bmatrix}_i \in \mathbb{R}^{2 + 2 N_C}$. The VT flash contains one less equality constraint than the UV flash. Consequently, there is also one less Lagrange multiplier, i.e. $z_i \in \mathbb{R}^{2 + N_C}$. The manipulated inputs remain unchanged. However, there are no disturbance variables in the isothermal model because the temperature of the injected water, $T^\text{inj}$, is constant. Consequently, the isothermal and compositional model is also in the semi-explicit differential-algebraic form \eqref{eq:sedae}.

%% file: tex/opt.tex
\section{Dynamic optimization}\label{sec:opt}
We consider the dynamic optimization problem

\begin{subequations}\label{eq:ocp}
	\begin{align}
	\label{eq:ocp:obj}
	\min_{[x(t);y(t);z(t)]_{t_0}^{t_f},\{u_k\}_{k=0}^{N-1}} \quad
	\phi &= \int_{t_0}^{t_f} \Phi(y(t), u(t), d(t)) dt,
	\end{align}
	subject to
	\begin{align}
	\label{eq:ocp:ic}
	&x(t_0) = \hat{x}_0, \\
	\label{eq:ocp:alg}  
	&G(x(t), y(t), z(t)) = 0,  & t &\in [t_0, t_f], \\
	\label{eq:ocp:dyn}
	&\dot{x}(t) = F(y(t),u(t),d(t)),  \enskip & t &\in [t_0, t_f], \\
	\label{eq:ocp:inp}
	&u(t) = u_k, \, & t &\in [t_k,\, t_{k+1}[, \, & k &= 0, \ldots, N-1, \\
	\label{eq:ocp:dist}
	&d(t) = \hat{d}_k, \, & t &\in [t_k,\, t_{k+1}[, \, & k &= 0, \ldots, N-1, \\
	\label{eq:ocp:con}
	&\{u_k\}_{k=0}^{N-1} \in \mathcal{U}.
	\end{align}
\end{subequations}
$[x(t); y(t); z(t)]_{t_0}^{t_f}$ is a vector of dependent decision variables, and $\{u_k\}_{k=0}^{N-1}$ are independent decision variables. $\hat{x}_0$ is an estimate of the initial states, and $\{\hat{d}_k\}_{k=0}^{N-1}$ are predictions of the disturbance variables. Both $\hat{x}_0$ and $\{\hat{d}_k\}_{k=0}^{N-1}$ are parameters in the optimization problem. $t_0$ is the initial time, and $t_N = t_f$ is the final time. $N$ is the number of control intervals.
\eqref{eq:ocp:ic} is a set of initial conditions for the semi-explicit differential-algebraic model equations \eqref{eq:ocp:alg}-\eqref{eq:ocp:dyn}. \eqref{eq:ocp:inp}-\eqref{eq:ocp:dist} are zero-order-hold (ZOH) parametrizations of the manipulated inputs and the disturbance variables. The constraints on the manipulated inputs \eqref{eq:ocp:con} are often bounds or linear constraints.

\subsection{The dynamic optimization algorithm}
We solve the dynamic optimization problem \eqref{eq:ocp} with the gradient-based algorithm described by Ritschel et al. \cite{Ritschel:etal:2017}. The algorithm is based on the single-shooting method which exploits that the initial value problem \eqref{eq:ocp:ic}-\eqref{eq:ocp:dyn}, subject to the ZOH parametrizations \eqref{eq:ocp:inp}-\eqref{eq:ocp:dist}, determines the dependent decision variables, $[x(t); y(t); z(t)]_{t_0}^{t_f}$, when $\{u_k\}_{k=0}^{N-1}$, $\hat{x}_0$, and $\{\hat{d}_k\}_{k=0}^{N-1}$ are specified. Consequently, the single-shooting method transcribes the infinite-dimensional dynamic optimization problem \eqref{eq:ocp} to a finite-dimensional numerical optimization problem in which the objective function requires the solution of the initial value problem (also referred to as a simulation). Efficient algorithms for solving numerical optimization problems require the gradients of the objective function. The dynamic optimization algorithm computes these gradients with the adjoint method. Furthermore, it uses Euler's implicit method to discretize the differential equations \eqref{eq:ocp:dyn}, and it solves the discretized differential equations and the algebraic equations in a simultaneous manner. The algorithm implements a simplified version of the time step selection scheme described by V\"{o}lcker et al. \cite{Volcker:etal:2010}, and it solves the involved linear systems with a block ILU(1) preconditioned GMRES method.

\subsection{Implementation}
We implement the dynamic optimization algorithm in C/C++. The implementation uses the open-source software DUNE for 1) grid management \cite{Bastian:etal:2008, Bastian:etal:2008b} and 2) solution of linear systems with the preconditioned GMRES method \cite{Blatt:Bastian:2007}. The involved thermodynamic functions (and their first and second order derivatives) are computed with C routines from ThermoLib. We use an SQLP algorithm \cite[Chap.~18]{Nocedal:Wright:2006}, from the commercial optimization software KNITRO 10.2, to solve the involved numerical optimization problem. Furthermore, we use C/C++ compilers from GCC.
In Section \ref{sec:ex}, we present performance tests which are carried out on a 64-bit workstation with 15.6~GB memory and four Intel Core i7 3.60~GHz cores. The workstation uses the Ubuntu 16.04 operating system. Furthermore, it has a shared level 3 cache of 8192~KB, and each core has a 256~KB level 2 cache and a 64~KB level 1 cache.

%% file: tex/ex.tex
\section{Numerical examples}\label{sec:ex}
In this section, we present numerical examples of optimized thermal and isothermal waterflooding strategies computed with the dynamic optimization algorithm described in Section \ref{sec:opt}. Furthermore, we discuss the computational performance of the algorithm in terms of various key performance indicators (KPIs).

\subsection{Optimized waterflooding strategies}\label{sec:ex:strat}
\begin{figure}
	\centering
	\includegraphics[width=0.45\textwidth]{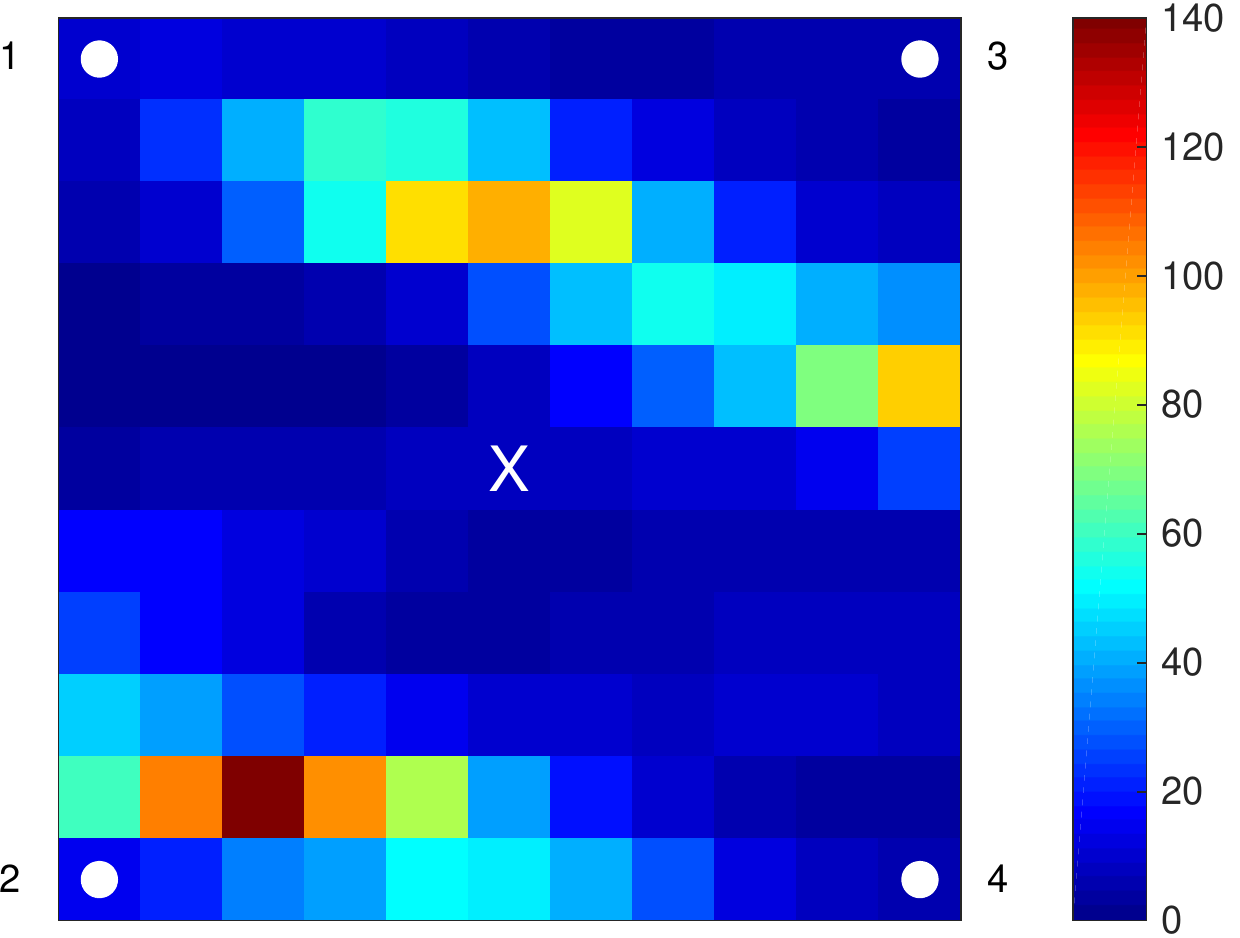}
	\caption{The permeability field [mD] and the locations of the injection and production wells. The white circles indicate the locations of the injectors, and the white X indicates the location of the producer.}
	\label{fig:perm}
\end{figure}

We consider a $110\times110\times10$~m reservoir which is initially at 50$^\circ$C. The oil and the gas phases consist of methane, ethane, propane, n-heptane, and hydrogen sulfide. We discretize the reservoir with an $11\times11\times1$ grid. The objective in the dynamic optimization of the thermal and isothermal waterflooding strategies is to maximize the total oil production over a three-year period. The decision variables are the BHPs of four injection wells and a single production well. There are 12 control intervals per year which results in a total of 36 control intervals, i.e. 36 decision variables per well. The locations of the wells are shown in Fig. \ref{fig:perm} together with the heterogeneous (and isotropic) permeability field. For simplicity, we assume that there is no heat loss to the surroundings, that the rock is incompressible, and that the porosity field is homogeneous. The porosity is $0.25$. The thermal heat capacity of the rock is $0.92$~$\text{kJ}/(\text{kg}\cdot\text{K})$, and the thermal conductivity is $2.5$~$\text{W}/(\text{m}\cdot\text{K})$ which resemble the properties of sandstone \cite[Chap.~2]{Eppelbaum:etal:2014}.
The BHPs of the injectors are constrained to the interval $[10~\text{MPa}, 12~\text{MPa}]$, and the BHP of the producer is constrained to the interval $[9~\text{MPa}, 10~\text{MPa}]$. The injected water is at $90^\circ$C in the thermal strategy and at $50^\circ$C in the isothermal strategy.

\begin{figure}
	\centering
	\includegraphics[width=0.95\textwidth]{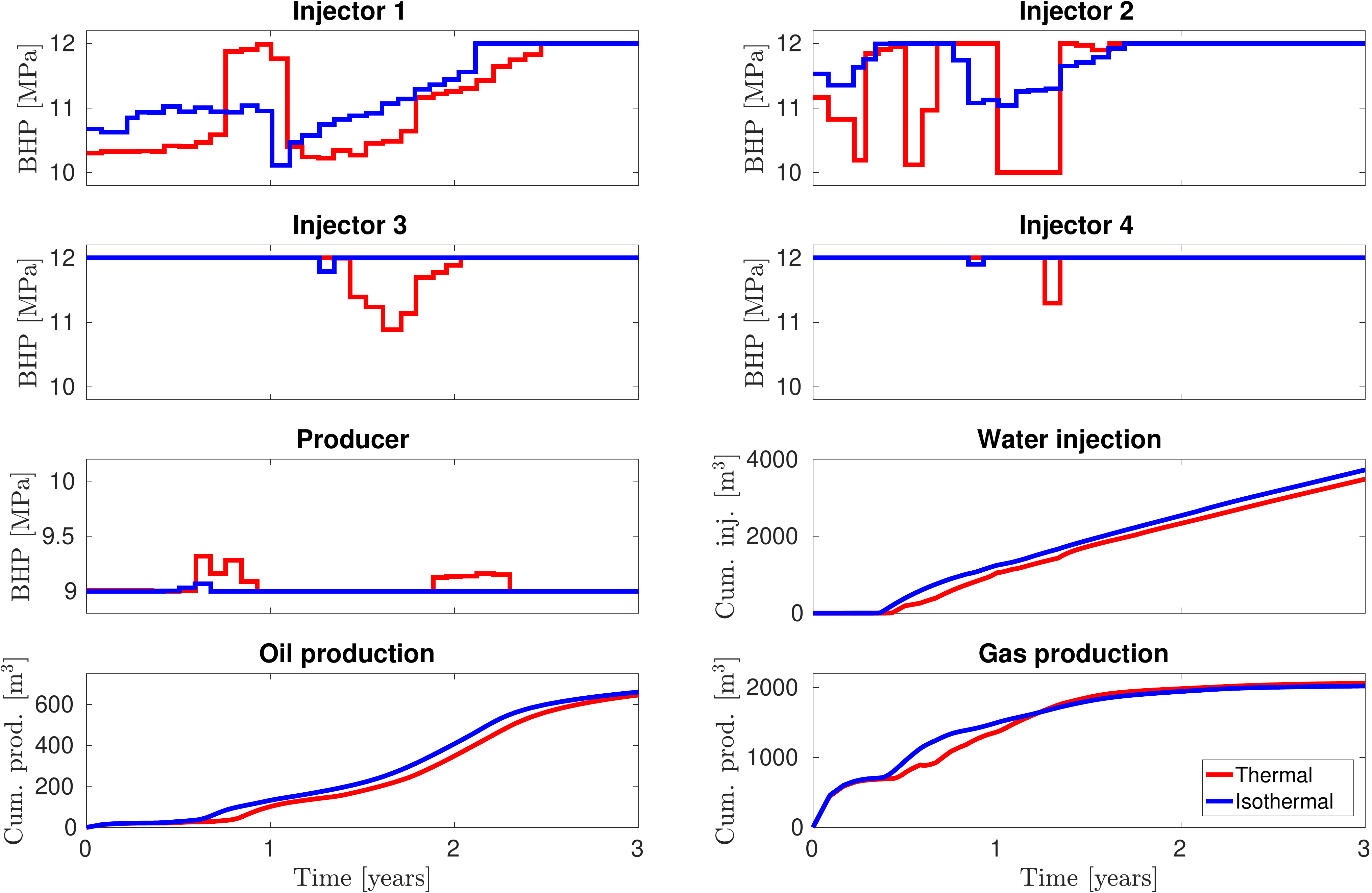}
	\caption{Optimized BHPs of the injectors and the cumulative volumetric water injection, oil production, and gas production.}
	\label{fig:opt}
\end{figure}

\begin{figure}
	\centering
	\includegraphics[width=0.211\textwidth]{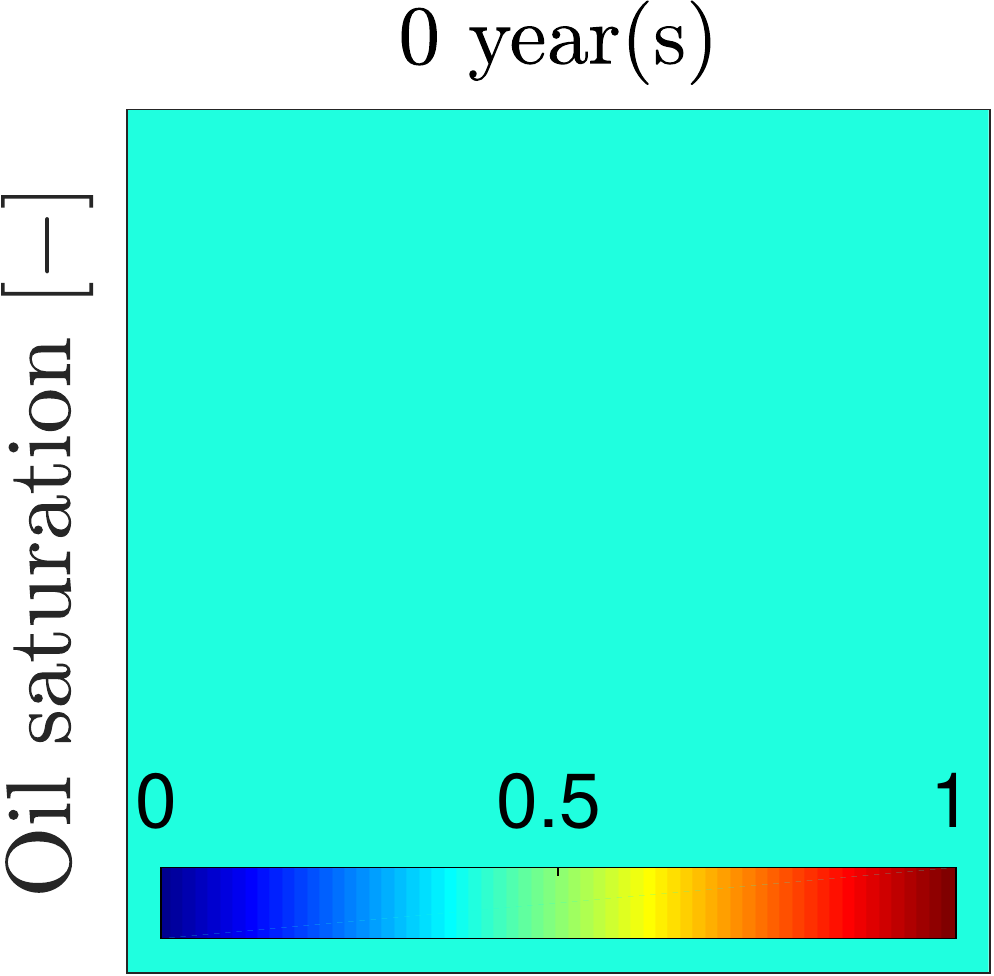}~
	\includegraphics[width=0.185\textwidth]{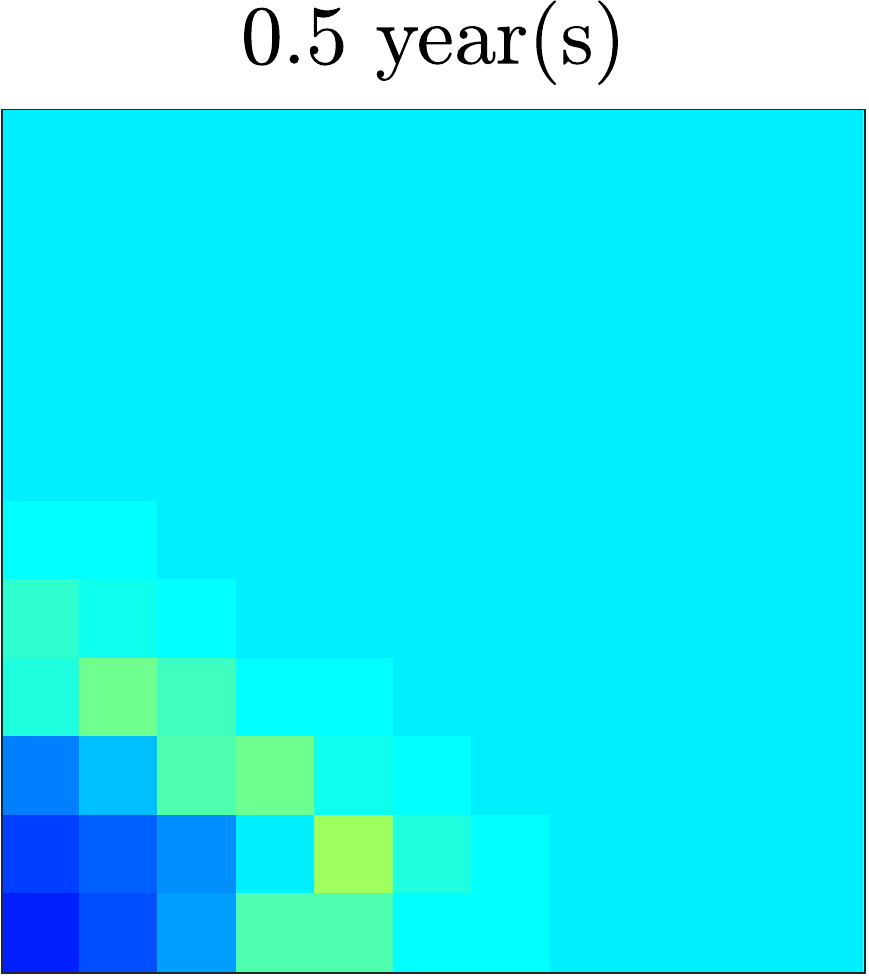}~
	\includegraphics[width=0.185\textwidth]{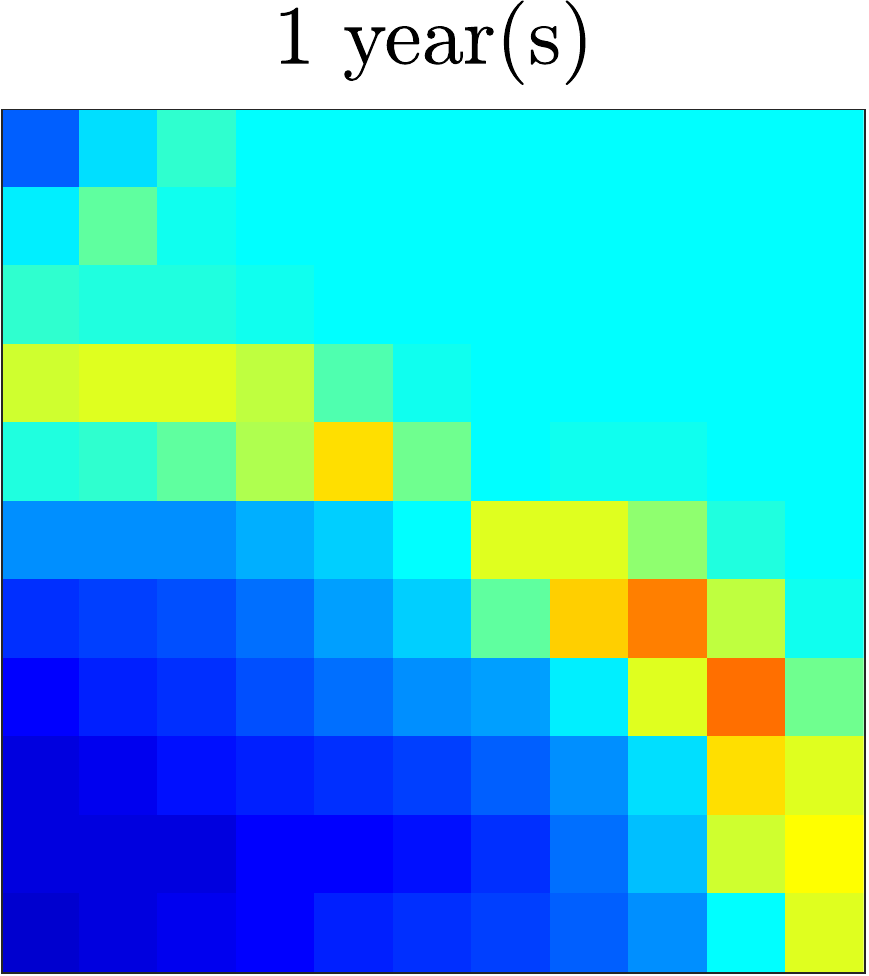}~
	\includegraphics[width=0.185\textwidth]{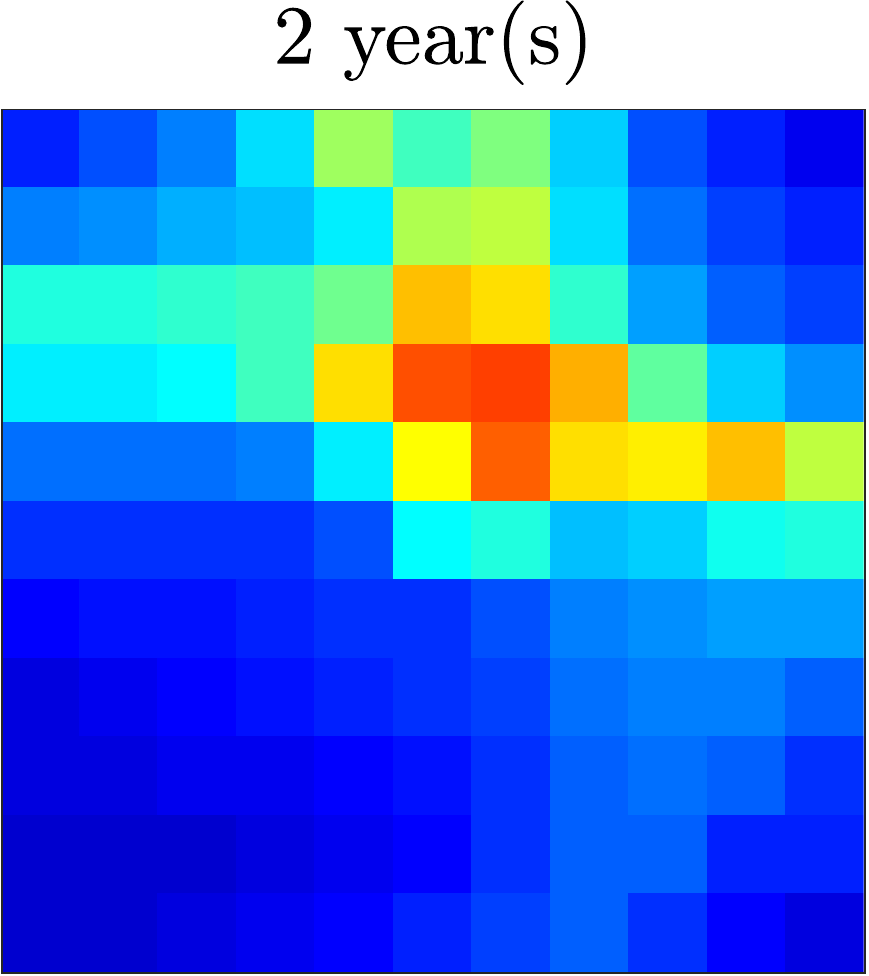}~
	\includegraphics[width=0.185\textwidth]{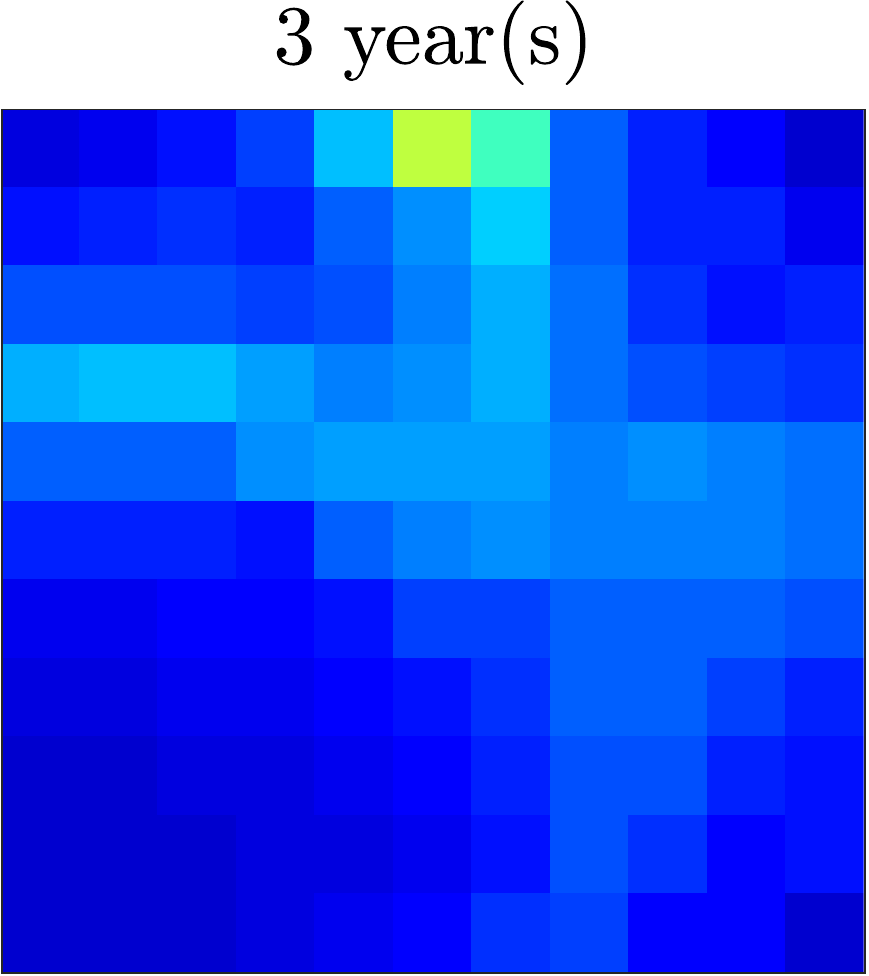} \\[2.5pt]
	\includegraphics[width=0.211\textwidth]{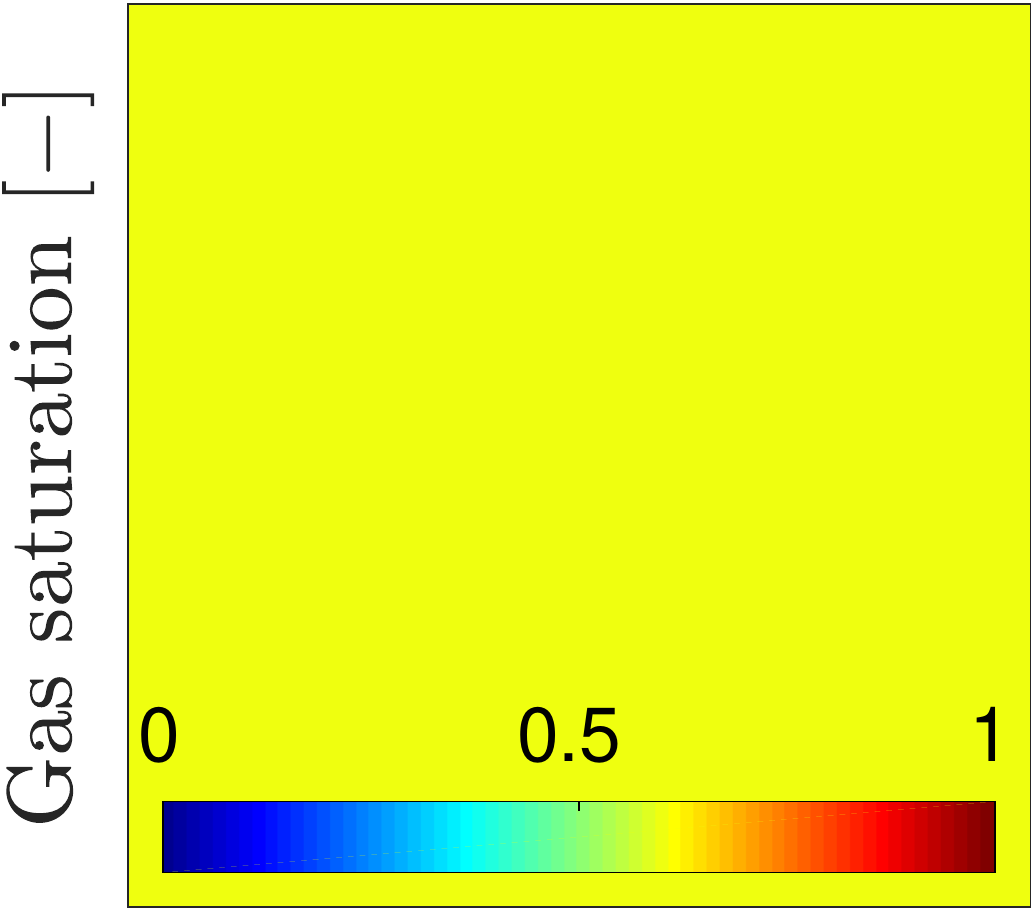}~
	\includegraphics[width=0.185\textwidth]{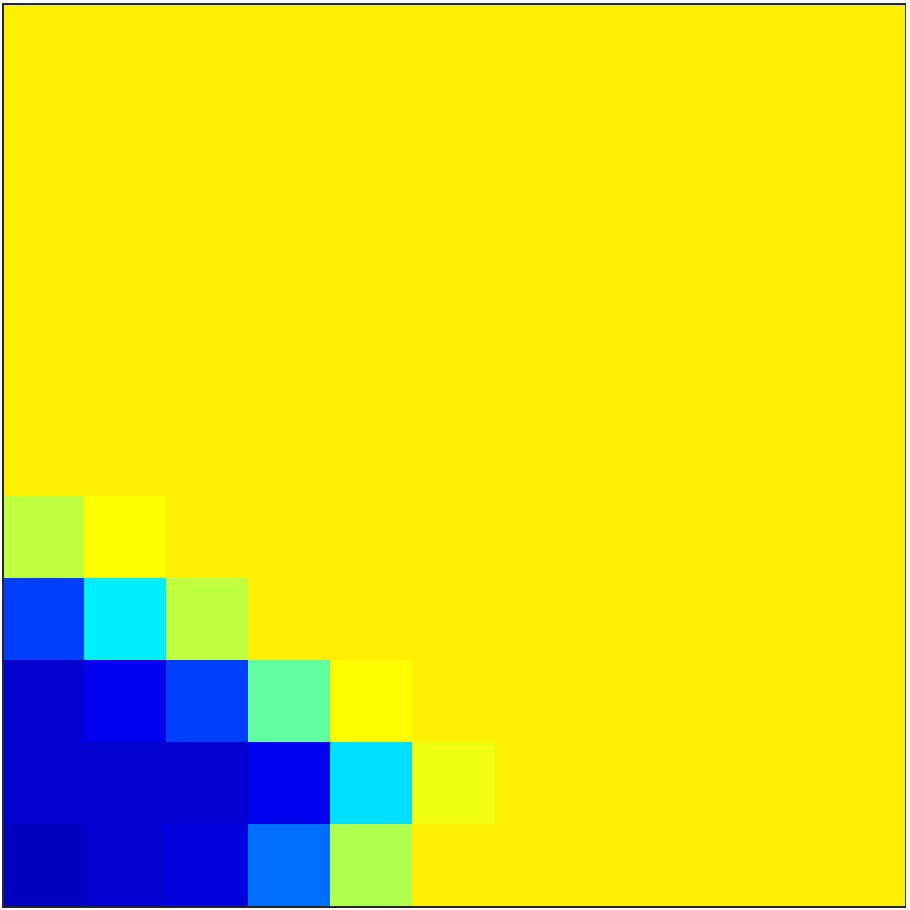}~
	\includegraphics[width=0.185\textwidth]{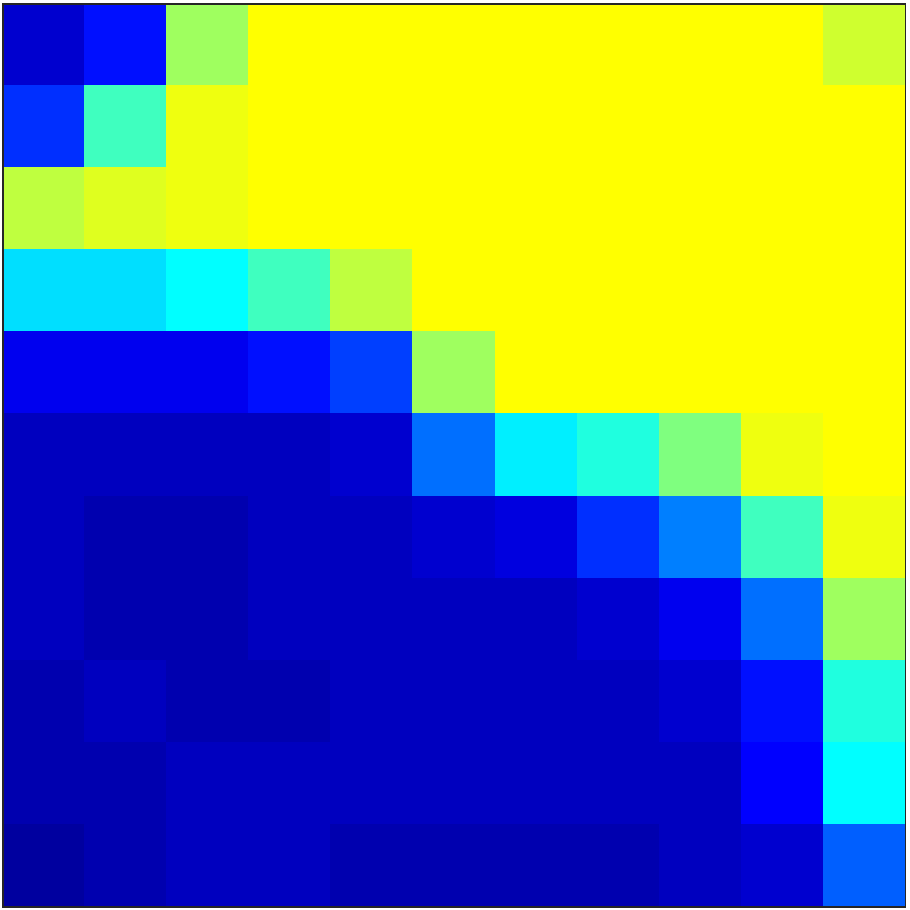}~
	\includegraphics[width=0.185\textwidth]{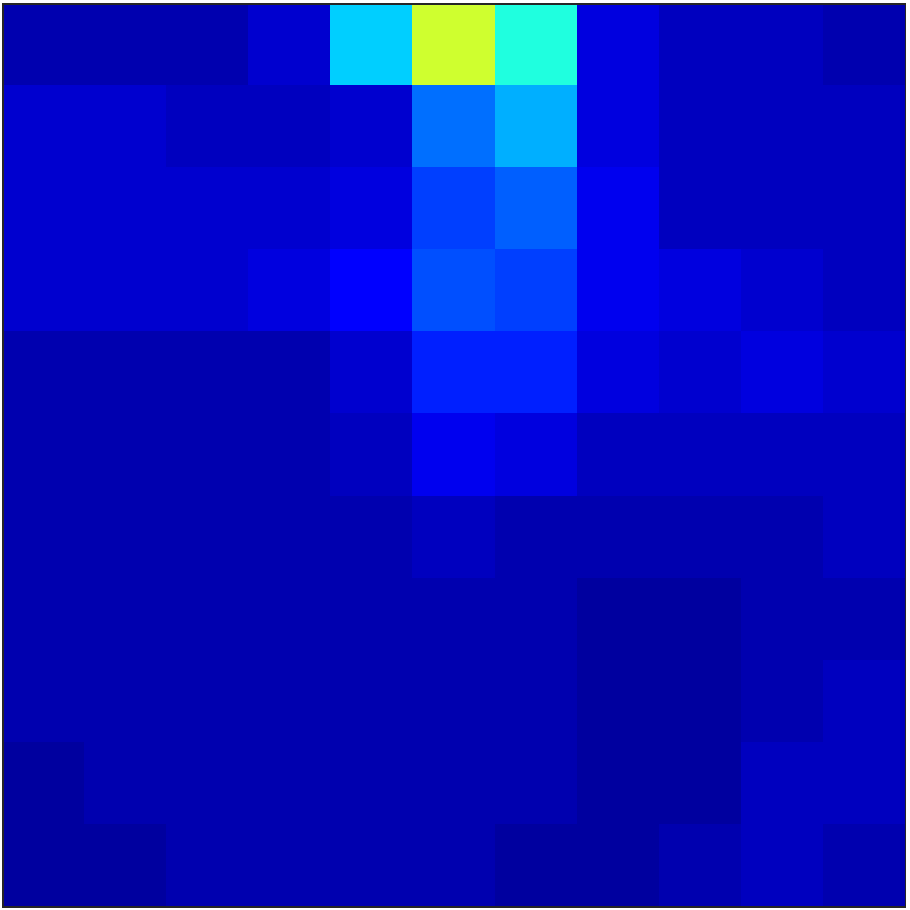}~
	\includegraphics[width=0.185\textwidth]{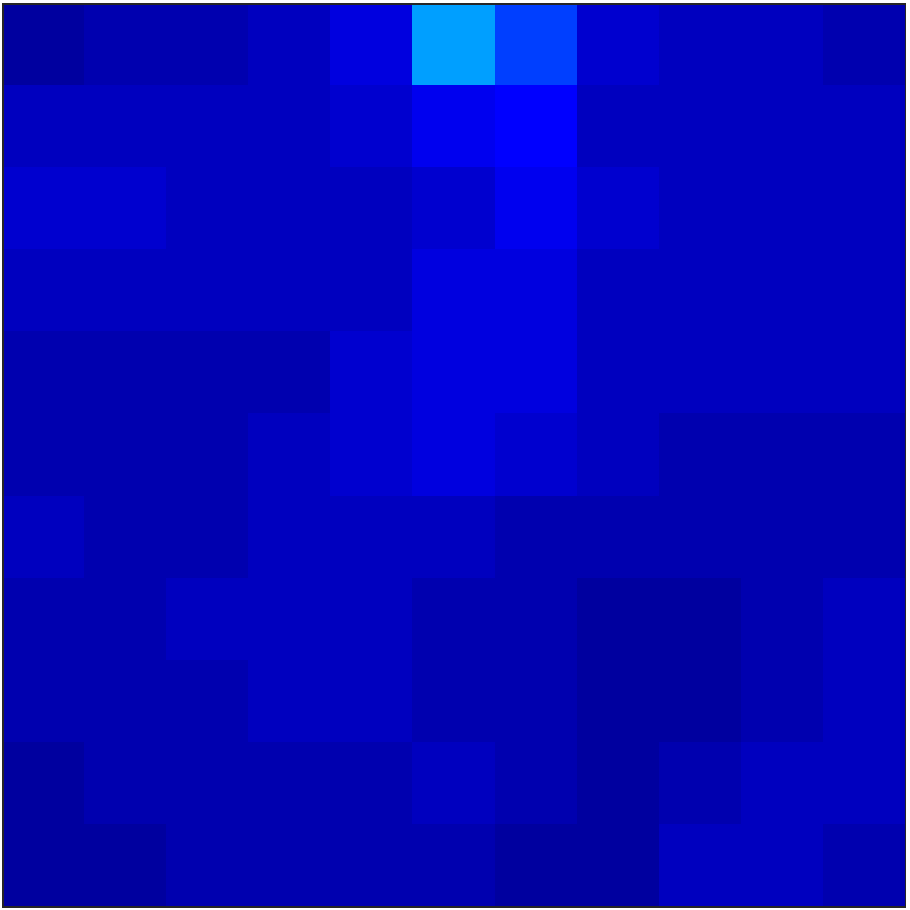} \\[2.5pt]
	\includegraphics[width=0.211\textwidth]{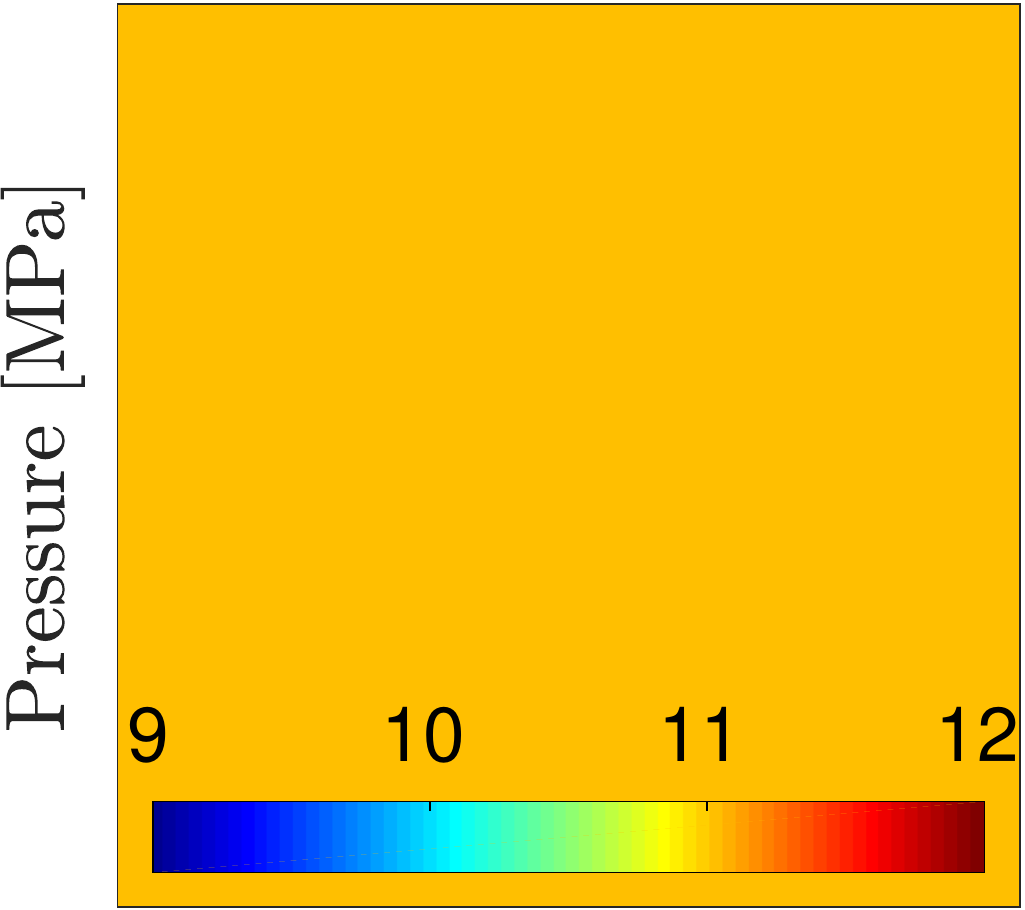}~
	\includegraphics[width=0.185\textwidth]{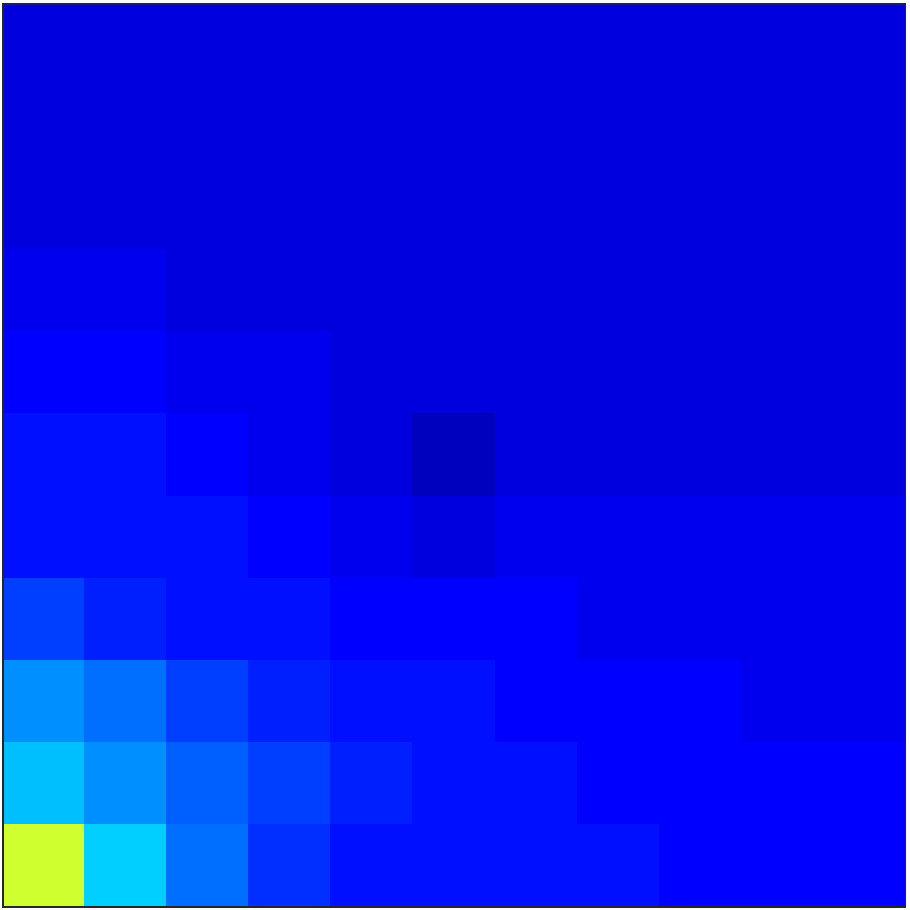}~
	\includegraphics[width=0.185\textwidth]{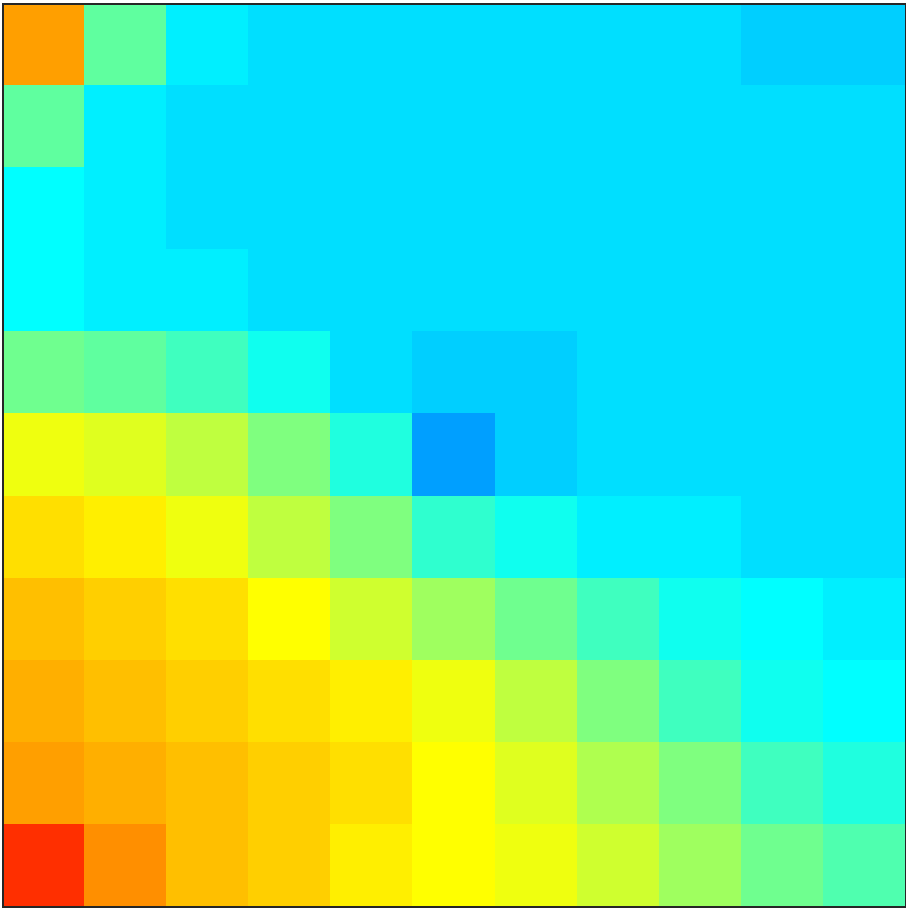}~
	\includegraphics[width=0.185\textwidth]{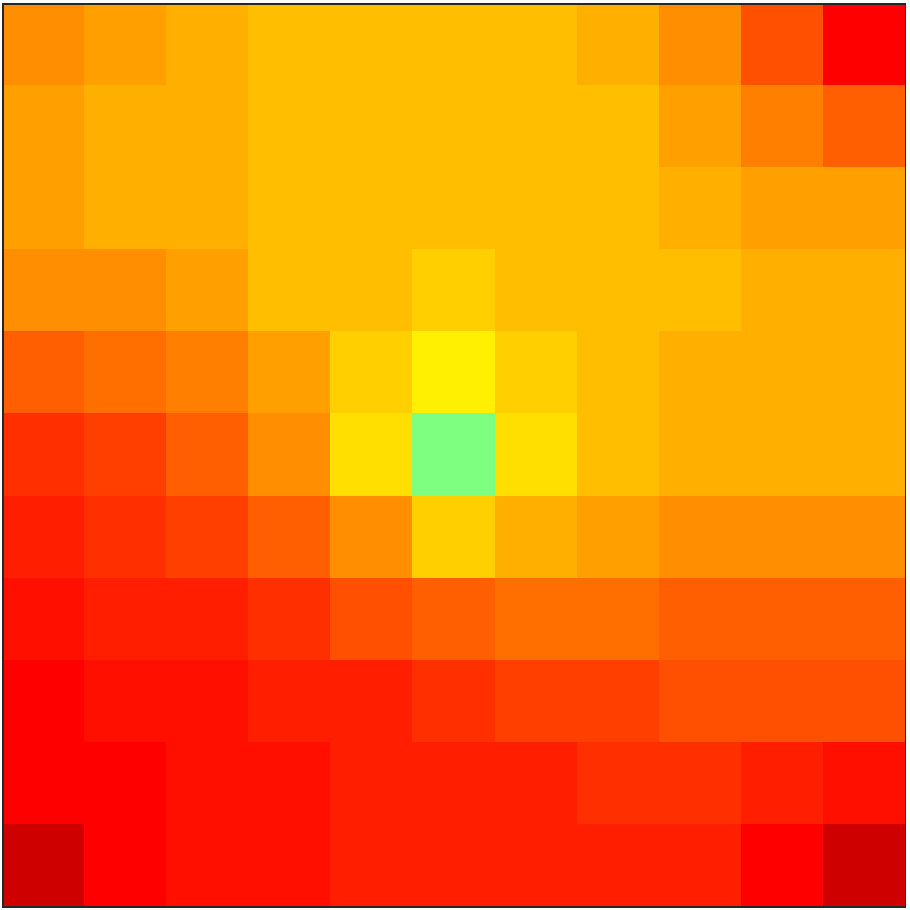}~
	\includegraphics[width=0.185\textwidth]{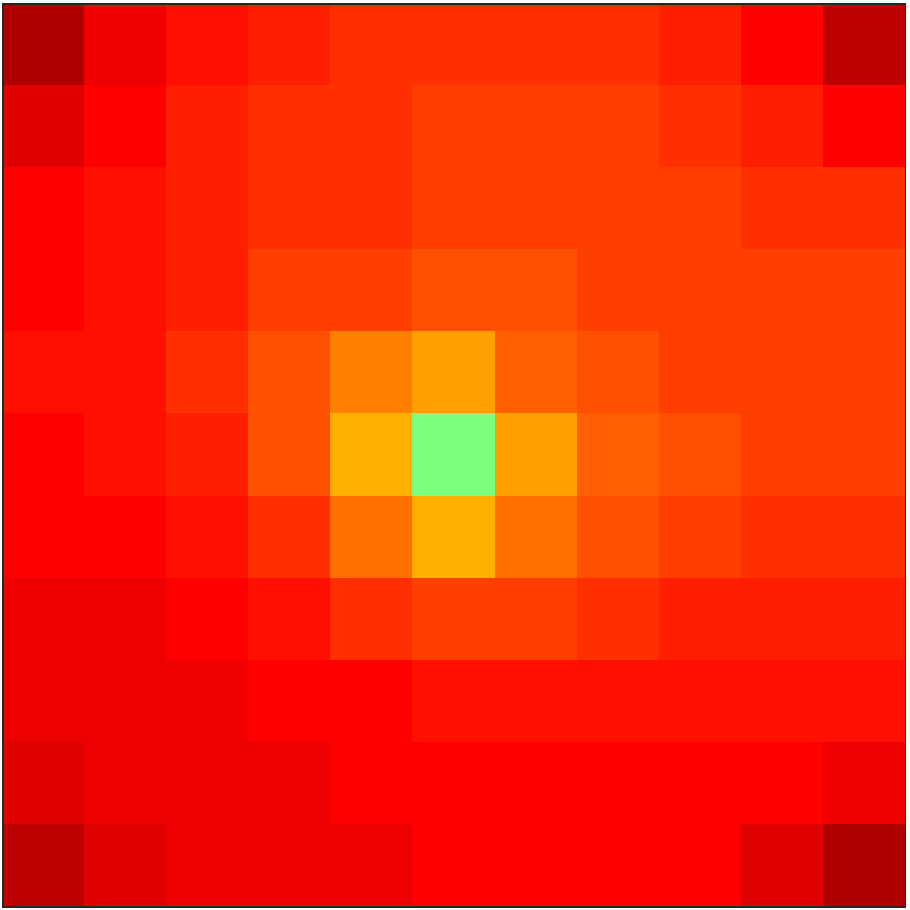} \\[2.5pt]
	\includegraphics[width=0.211\textwidth]{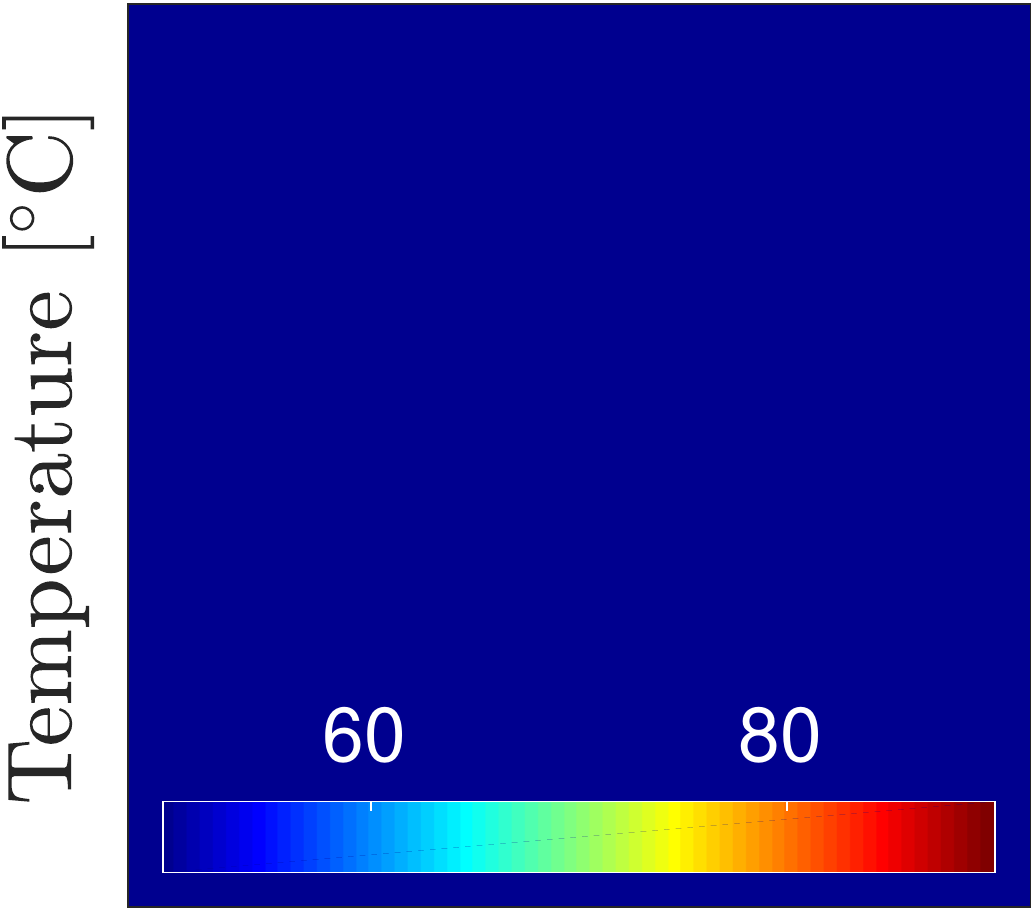}~
	\includegraphics[width=0.185\textwidth]{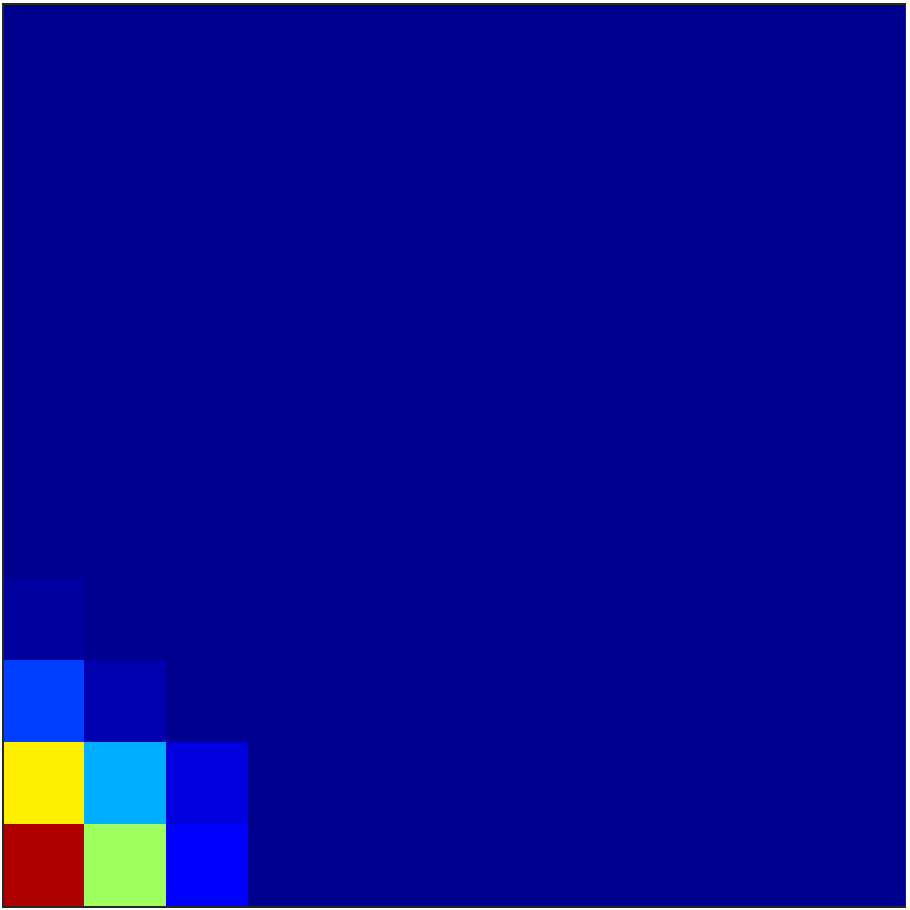}~
	\includegraphics[width=0.185\textwidth]{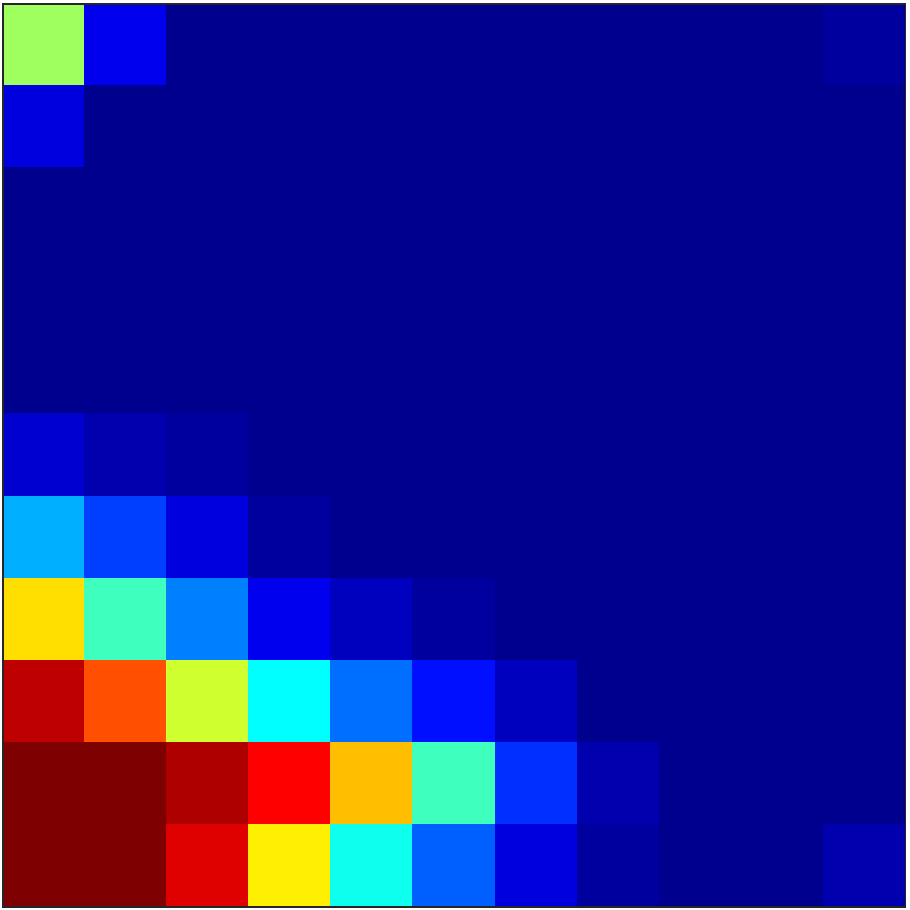}~
	\includegraphics[width=0.185\textwidth]{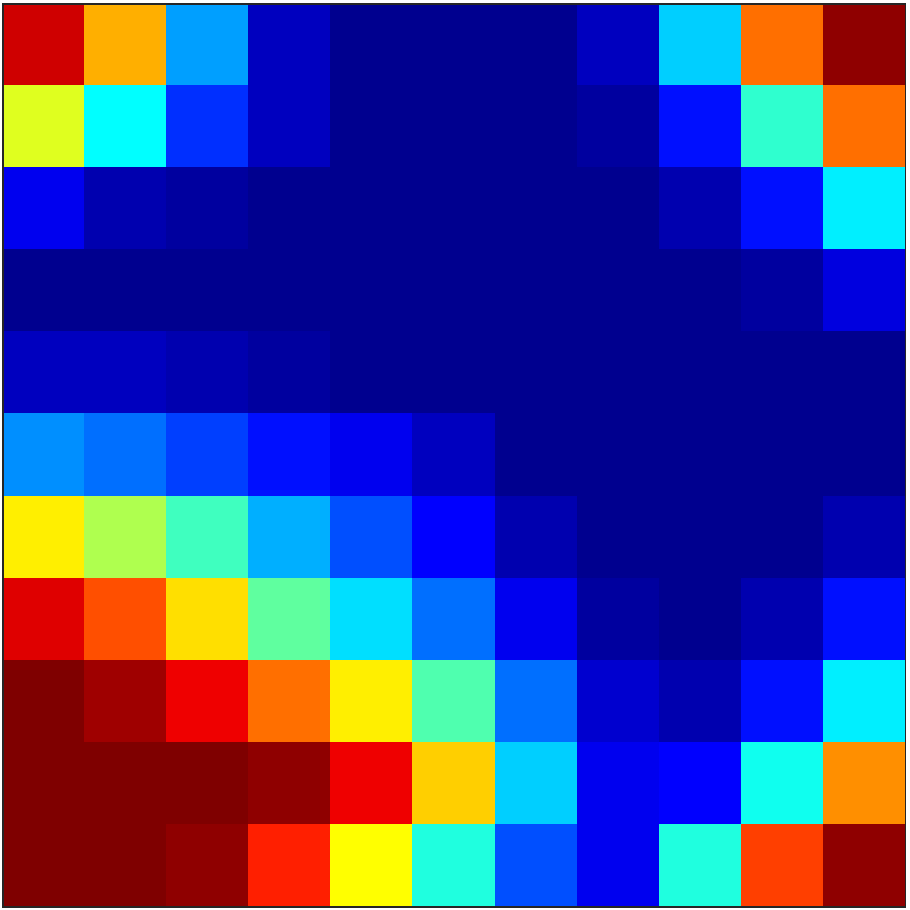}~
	\includegraphics[width=0.185\textwidth]{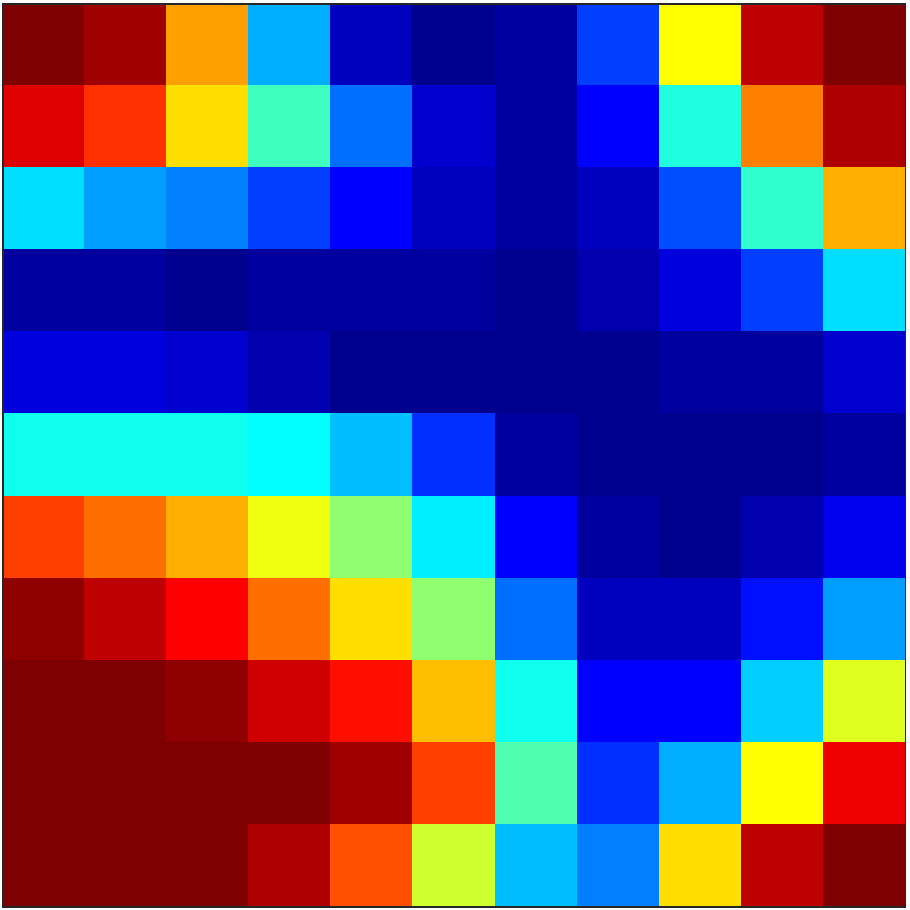}
	\caption{Optimized thermal waterflooding (90$^\circ$C water injected).}
	\label{fig:thermal}
\end{figure}

\begin{figure}
	\centering
	\includegraphics[width=0.211\textwidth]{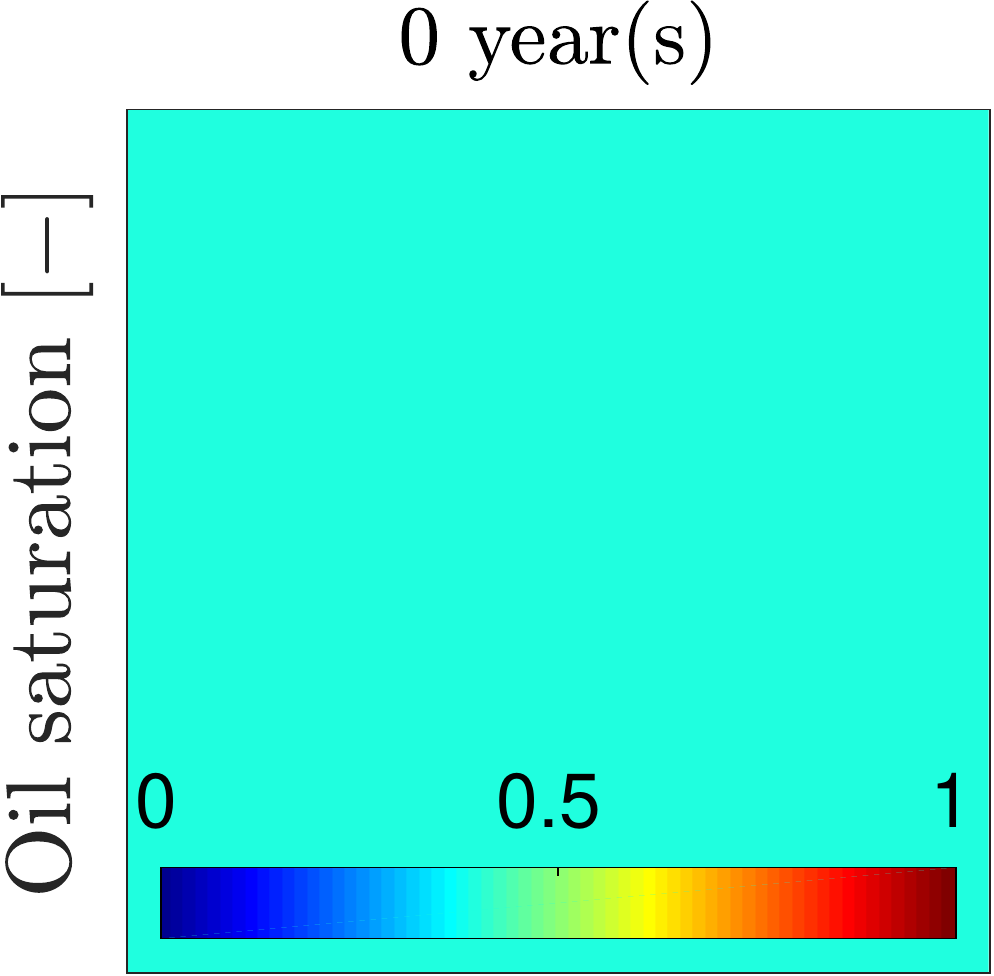}~
	\includegraphics[width=0.185\textwidth]{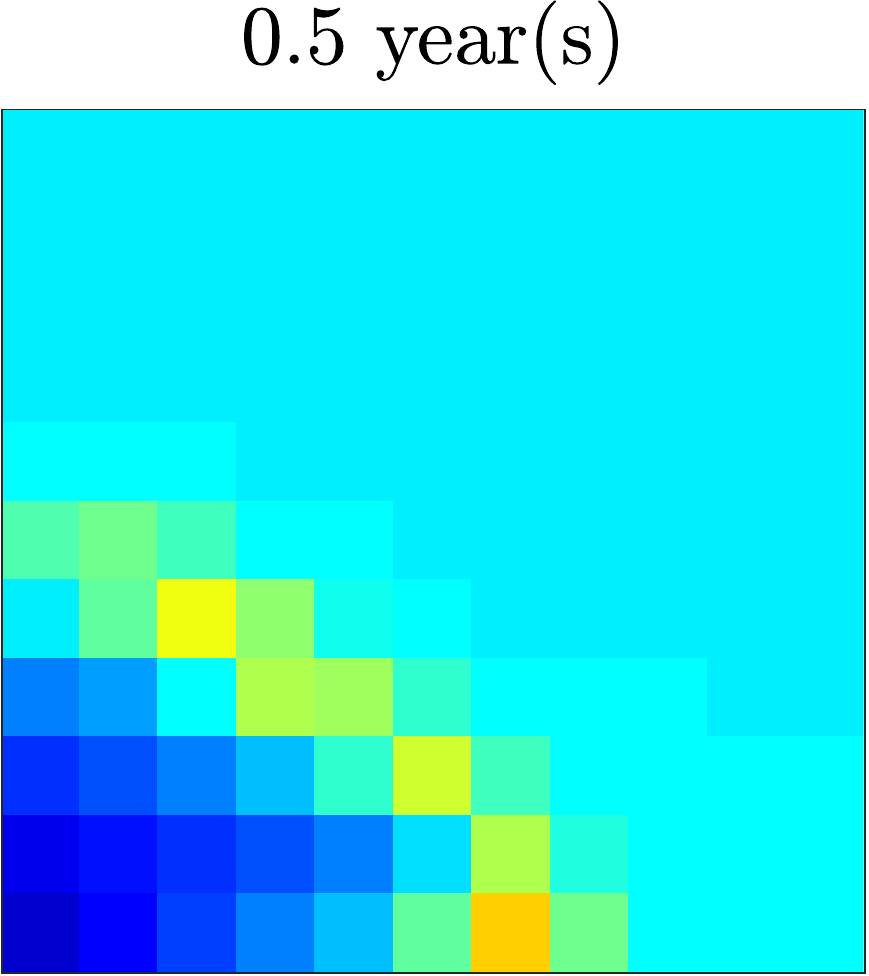}~
	\includegraphics[width=0.185\textwidth]{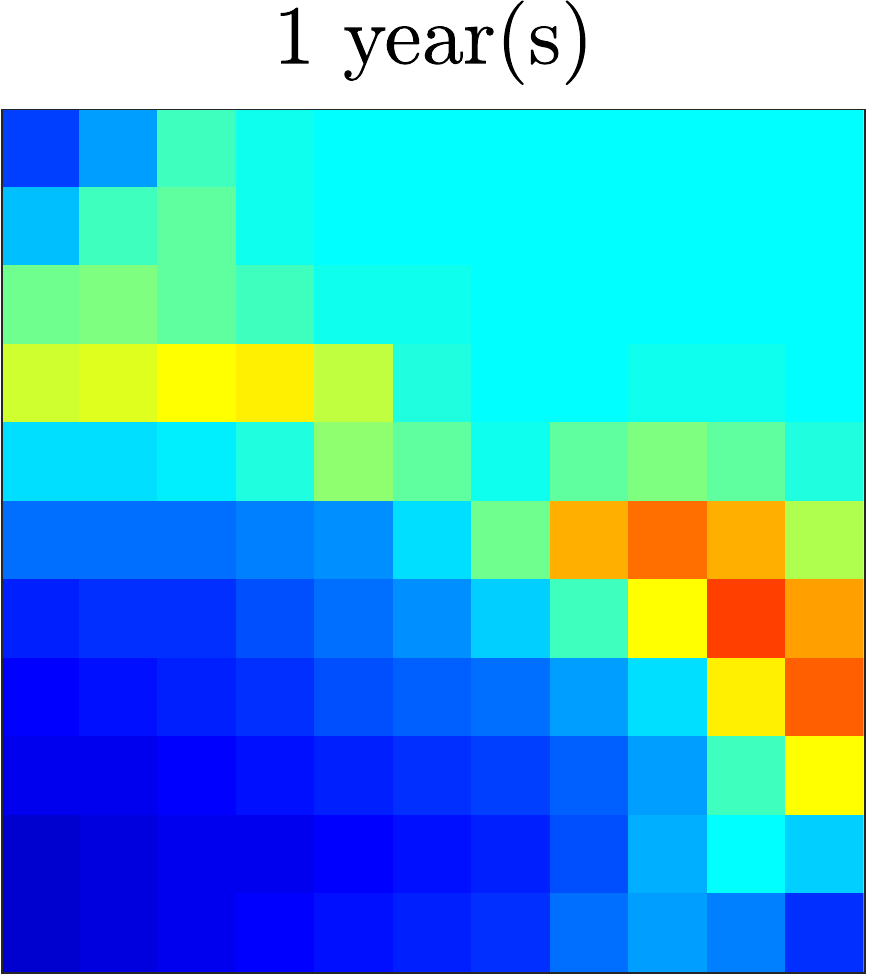}~
	\includegraphics[width=0.185\textwidth]{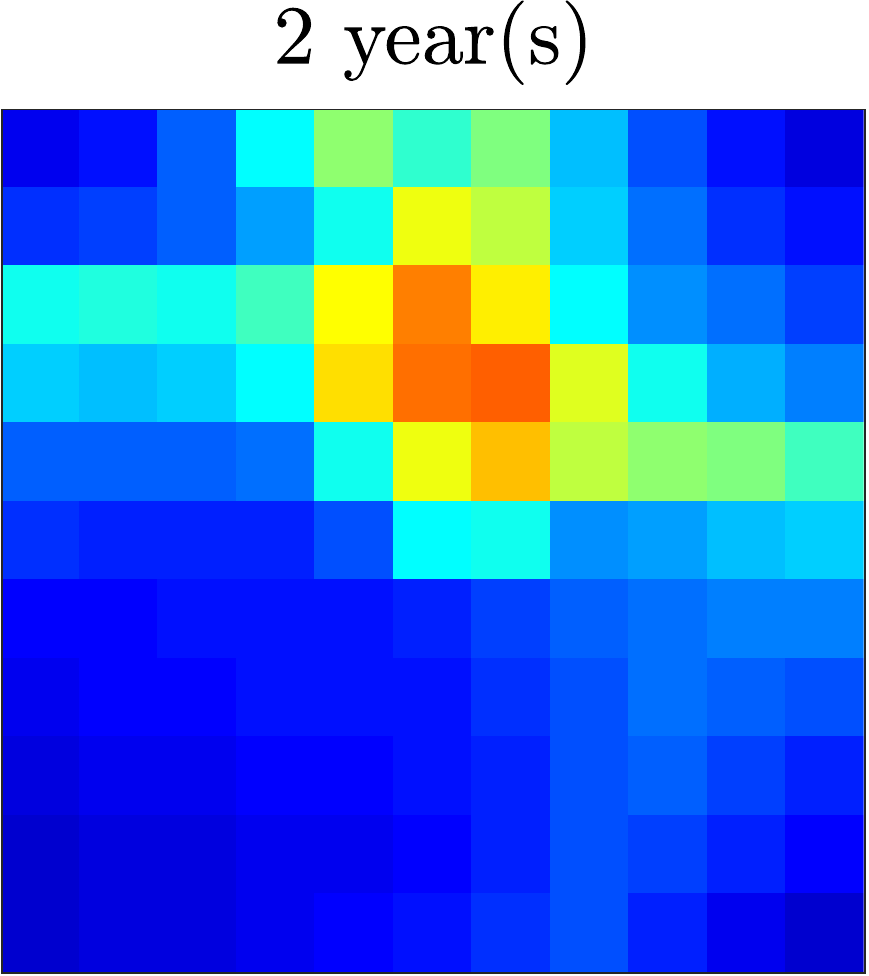}~
	\includegraphics[width=0.185\textwidth]{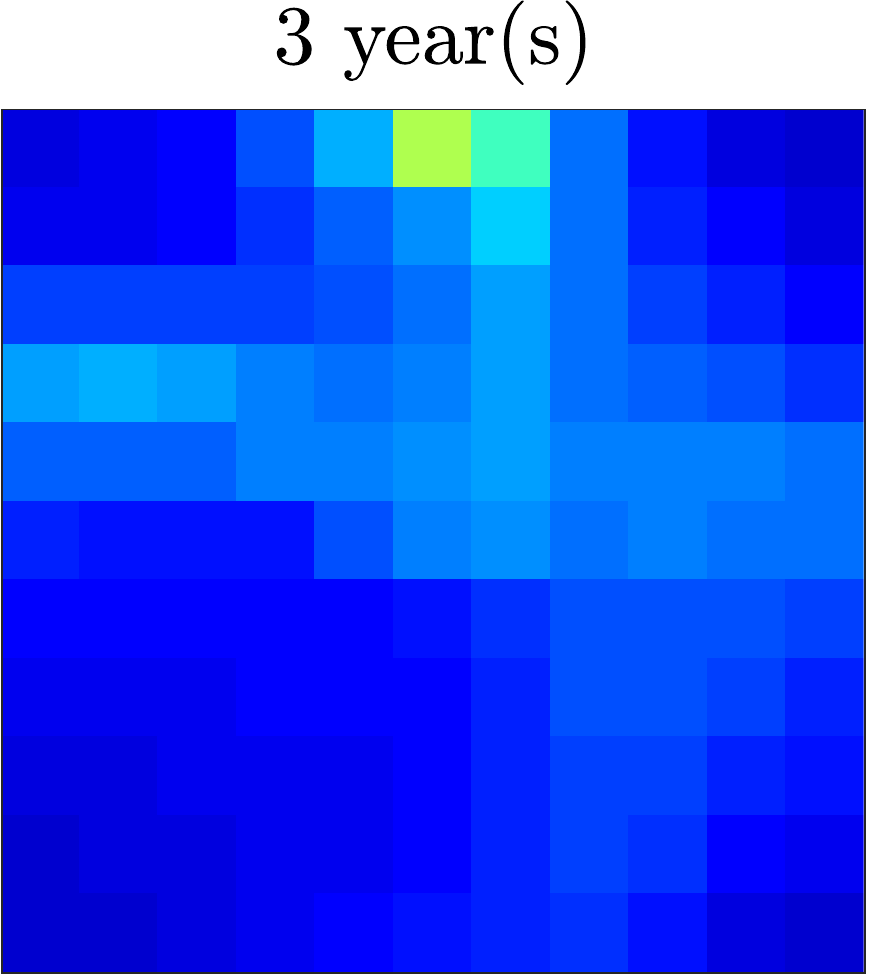} \\[2.5pt]
	\includegraphics[width=0.211\textwidth]{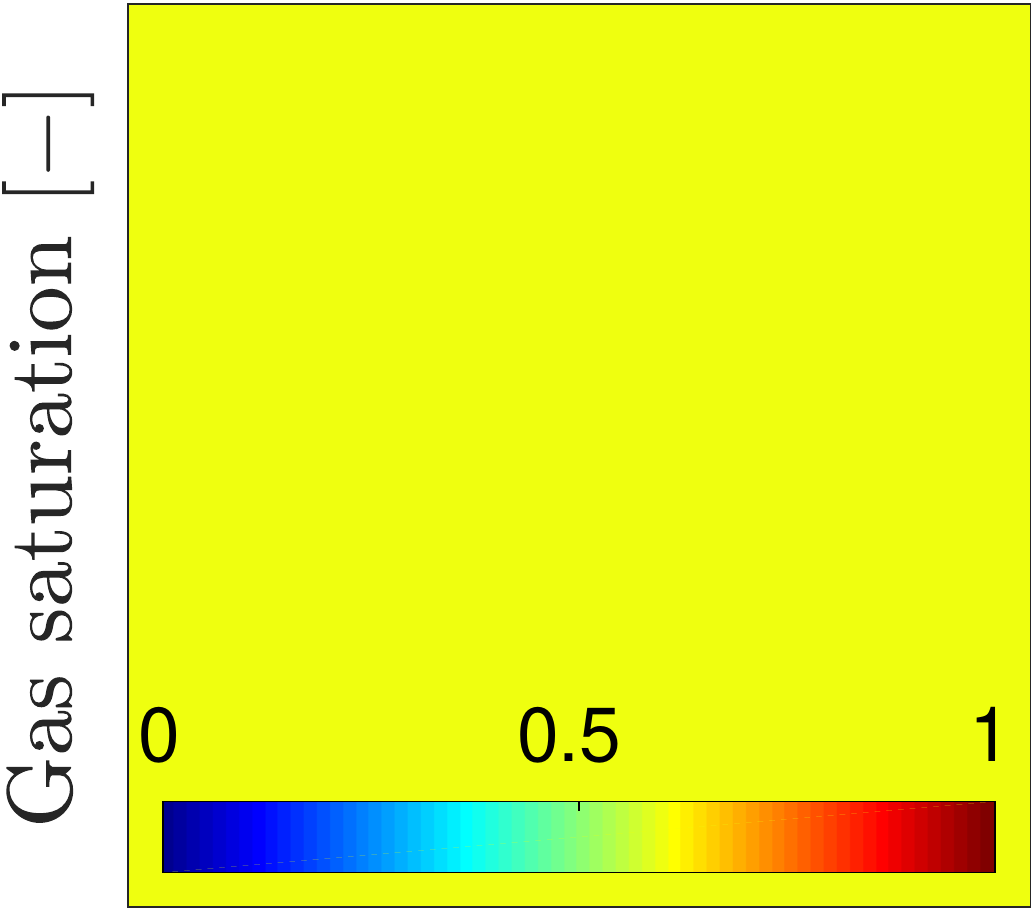}~
	\includegraphics[width=0.185\textwidth]{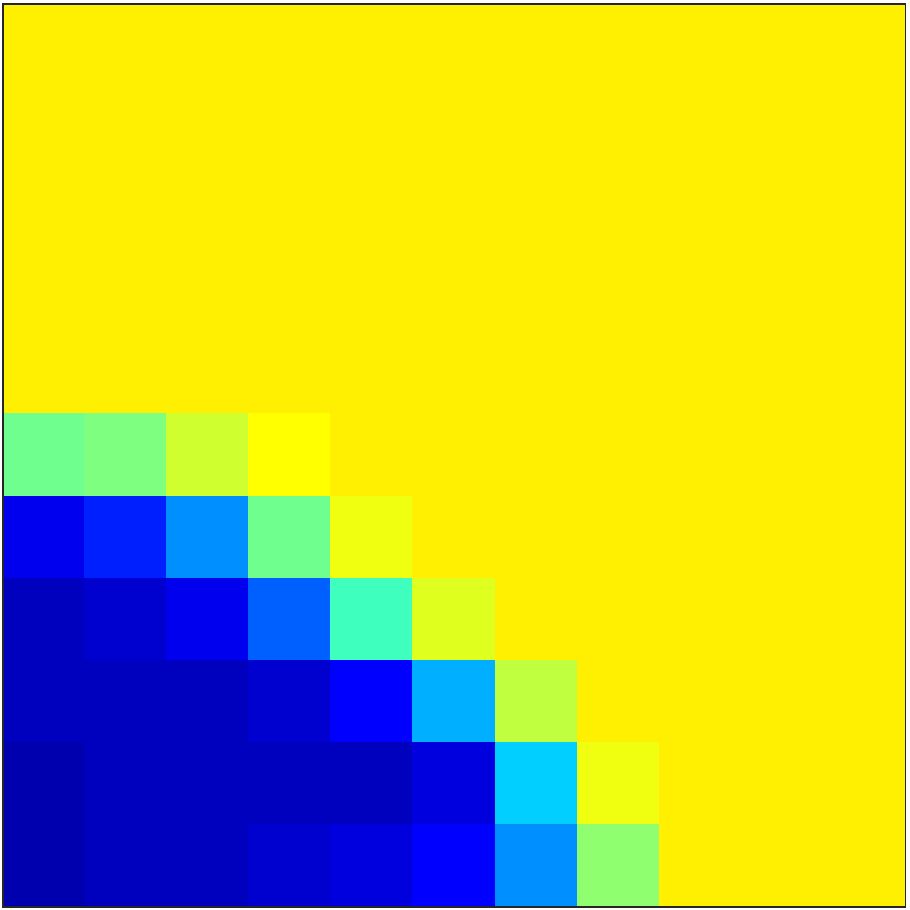}~
	\includegraphics[width=0.185\textwidth]{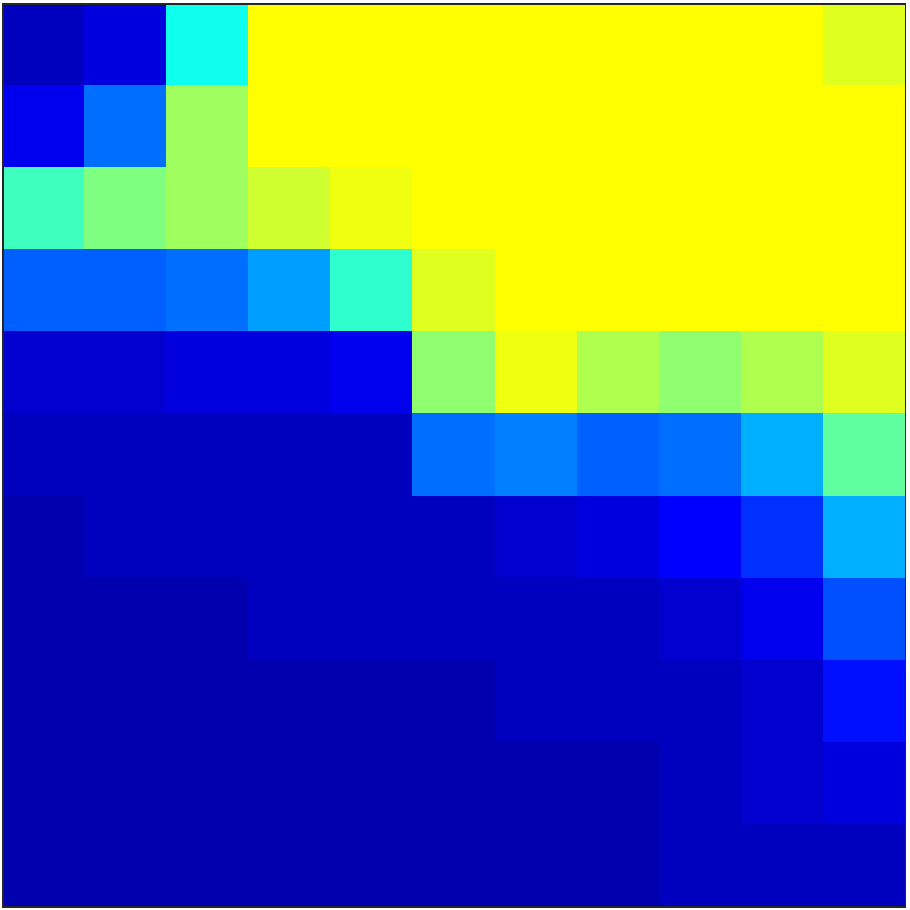}~
	\includegraphics[width=0.185\textwidth]{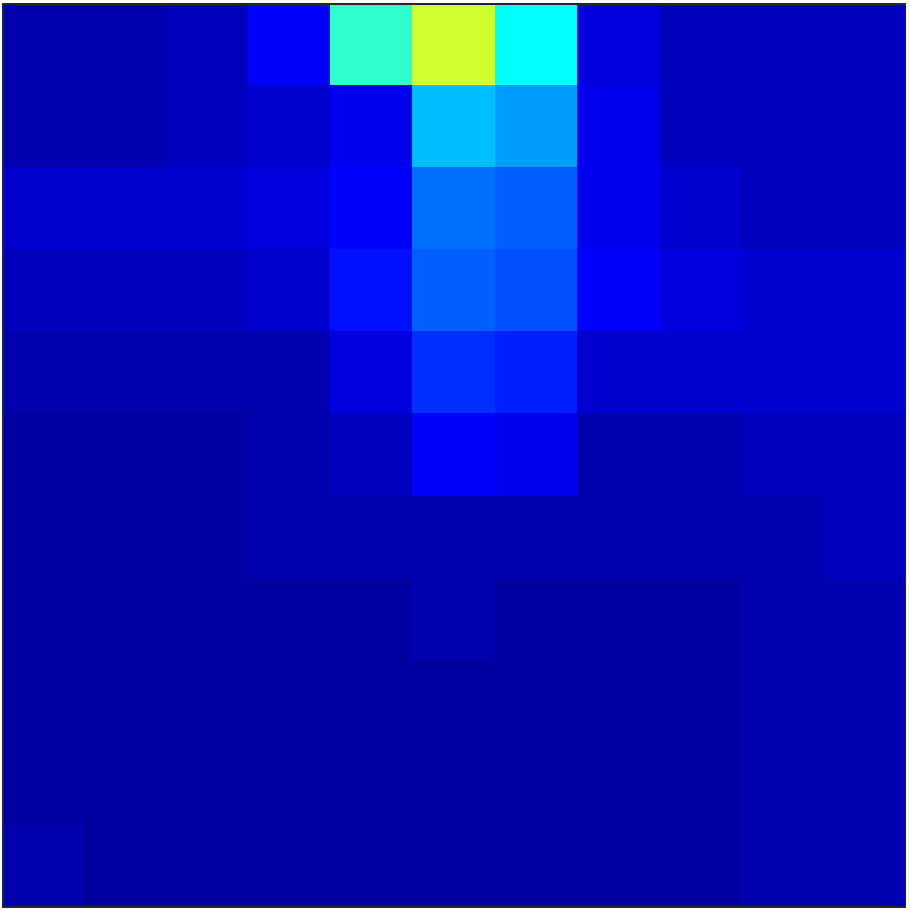}~
	\includegraphics[width=0.185\textwidth]{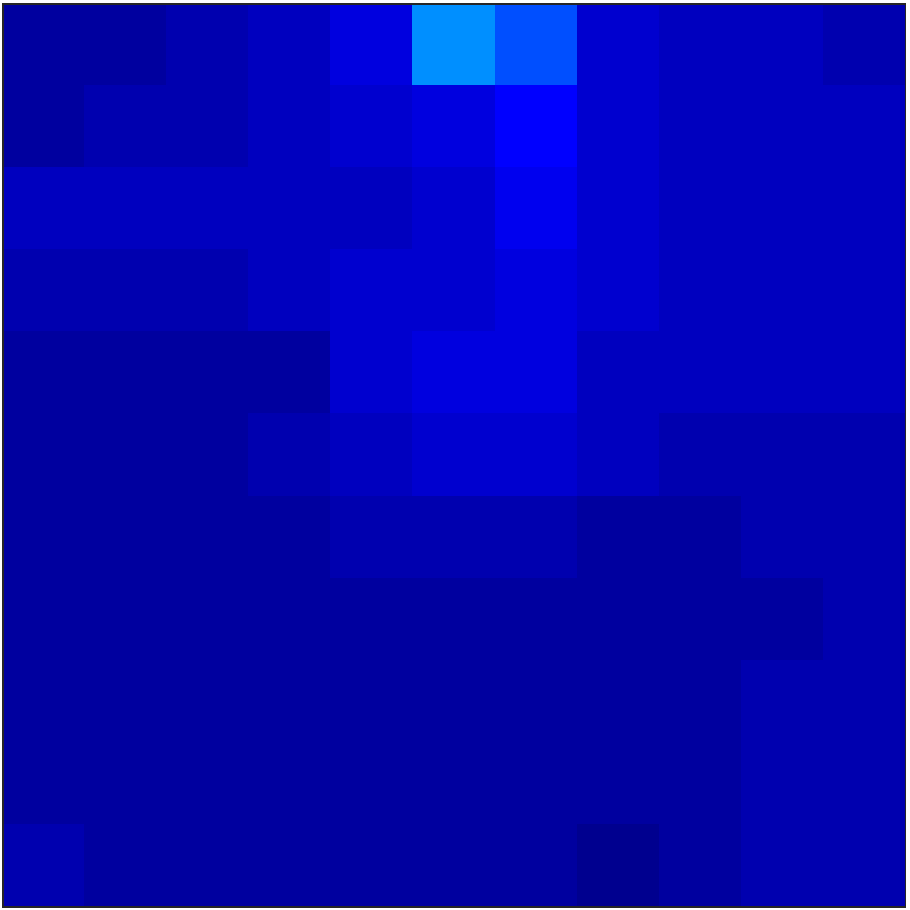} \\[2.5pt]
	\includegraphics[width=0.211\textwidth]{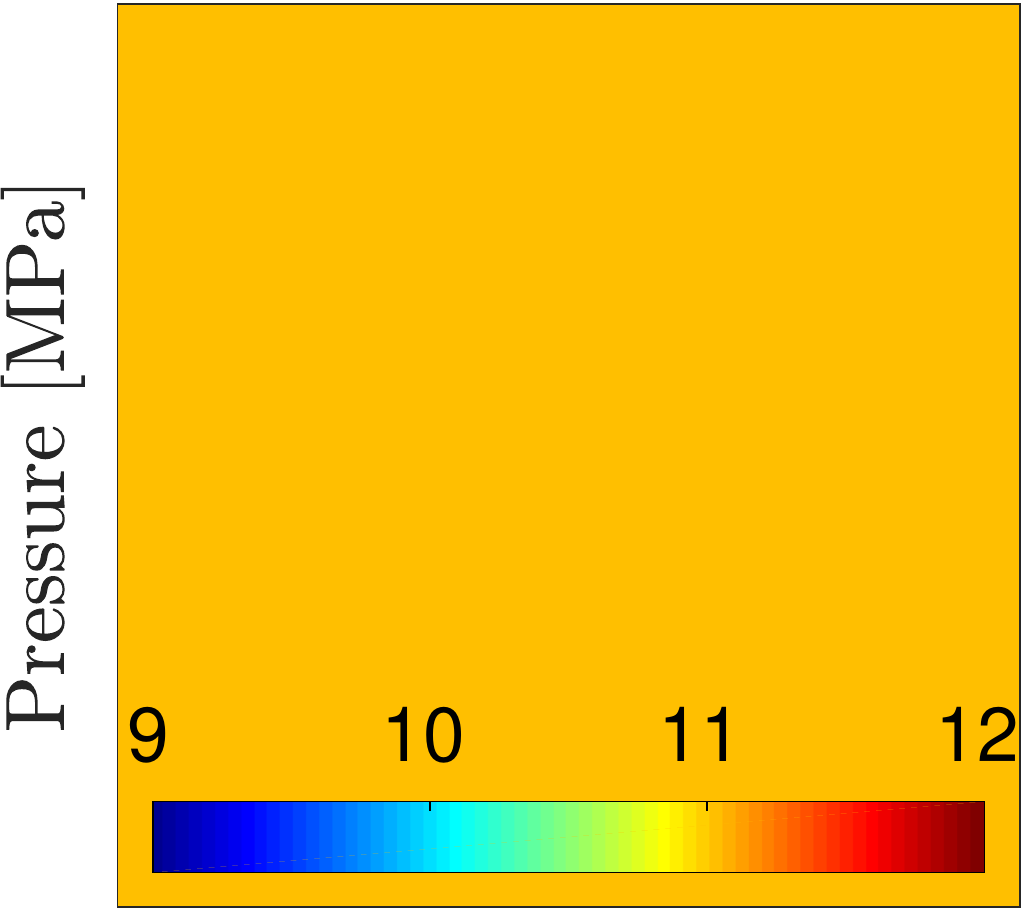}~
	\includegraphics[width=0.185\textwidth]{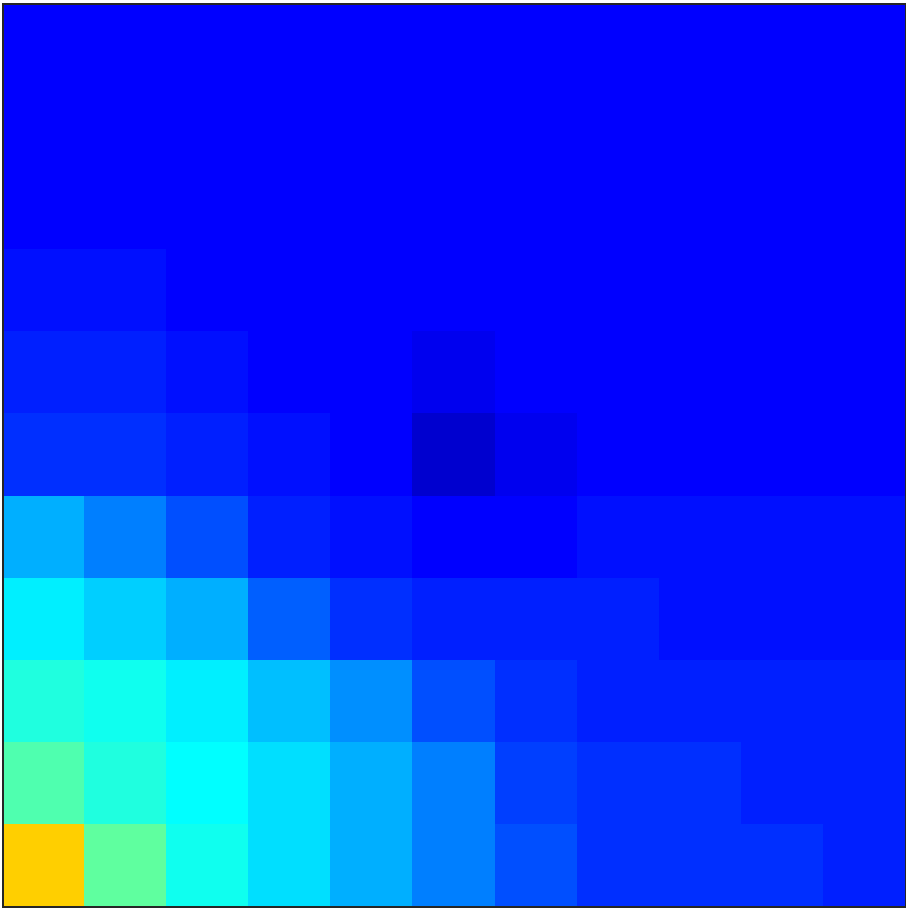}~
	\includegraphics[width=0.185\textwidth]{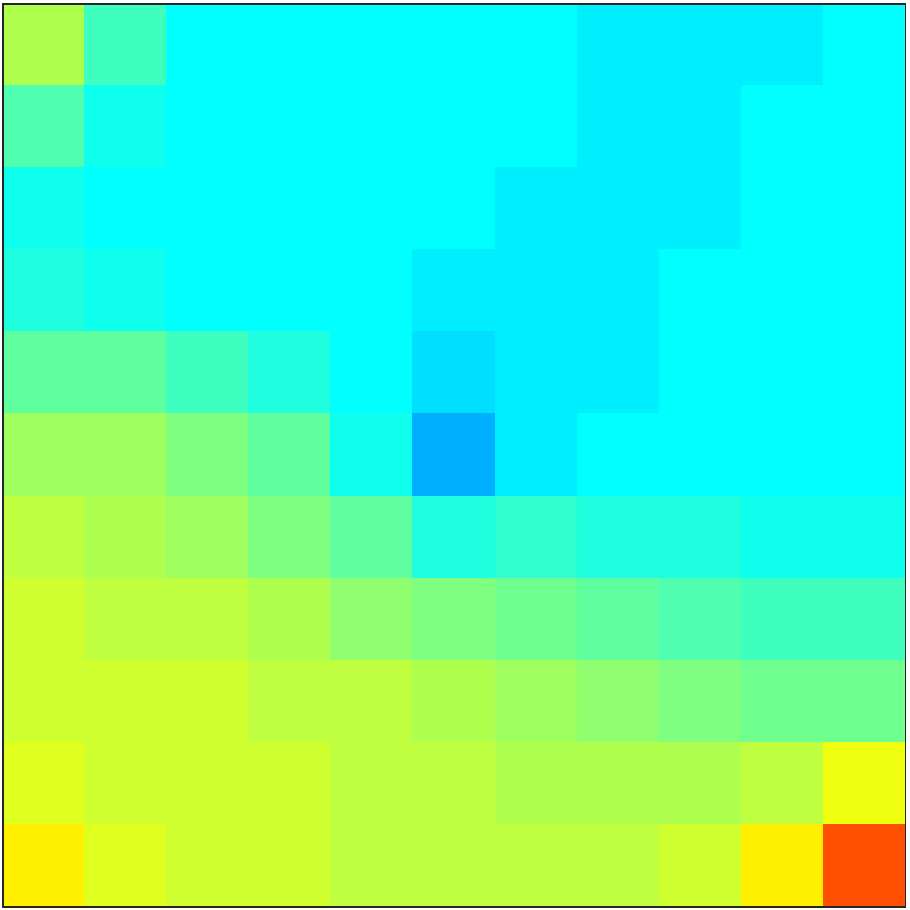}~
	\includegraphics[width=0.185\textwidth]{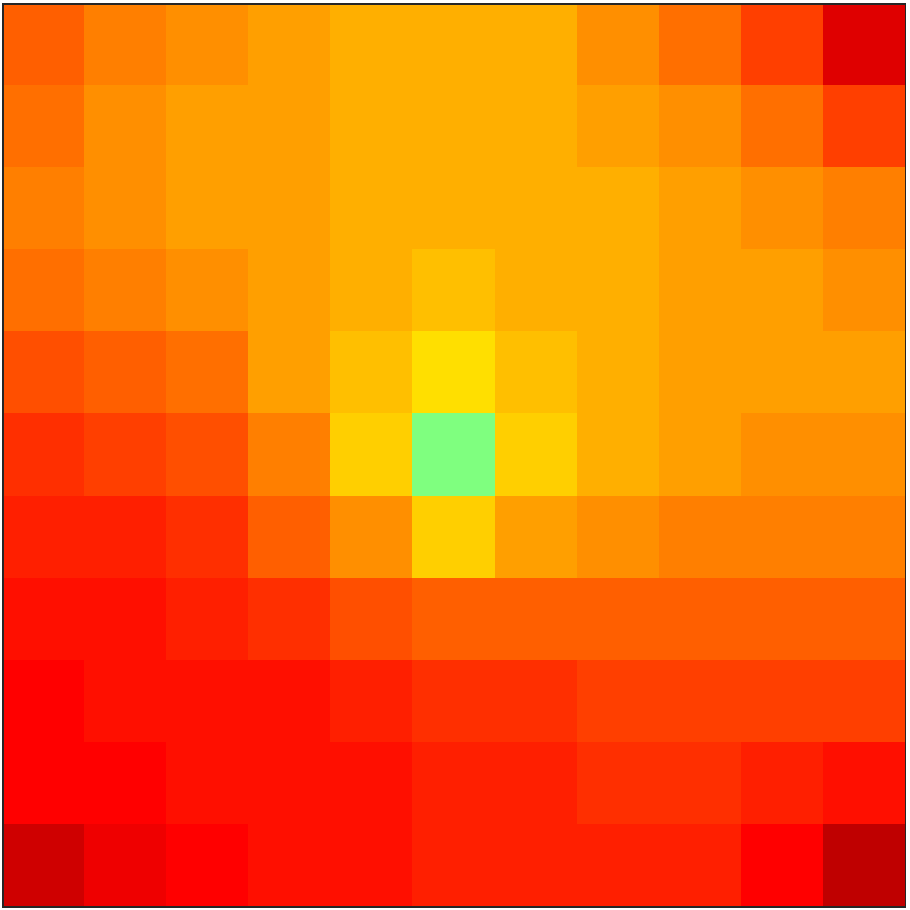}~
	\includegraphics[width=0.185\textwidth]{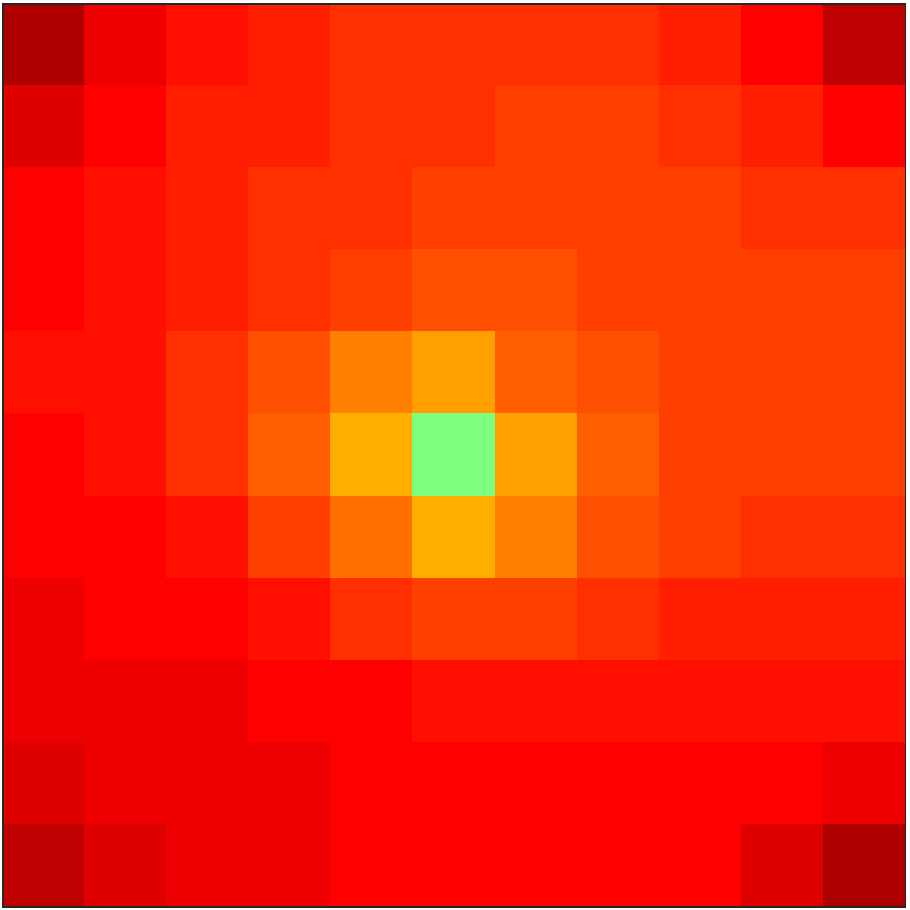}
	\caption{Optimized isothermal waterflooding (50$^\circ$C water injected).}
	\label{fig:isothermal}
\end{figure}

Fig. \ref{fig:opt} shows the injector and producer BHPs of the optimized thermal and isothermal strategies together with the cumulative volumetric injection of water, production of oil, and production of gas. Both strategies operate the producer close to the minimum BHP. Furthermore, they operate injector 3 and 4 close to the maximum BHP because they are located in very impermeable areas. The most significant differences between the two strategies are the BHPs of injector 1 and 2 which they both vary significantly. Compared to the isothermal strategy, the thermal strategy 1) injects slightly less water, 2) produces slightly less oil, and 3) produces slightly more gas.
Fig. \ref{fig:thermal} and \ref{fig:isothermal} illustrate the thermal and isothermal waterflooding processes in terms of the pressure and the oil and gas saturations in the reservoir. Fig. \ref{fig:thermal} also shows the temperature in the reservoir. The two figures suggest that it is challenging for the dynamic optimization algorithm to compute strategies that completely deplete the upper half of the reservoir during the three years of production.

In this example, we have considered an optimized thermal waterflooding strategy. However, thermal oil recovery strategies most often involve the injection of steam, e.g. steam-assisted gravity drainage (SAGD) is used to produce heavy oil in Canada and Venezuela \cite{Alvarado:Manrique:2010}. In such cases, the temperature of the injected steam would be a manipulated input together with the well BHP, and the objective function should include the cost of heating the steam.

\subsection{Key performance indicators}
Table \ref{tab:kpi} shows a number of key problem characteristics and KPIs for the dynamic optimization of the thermal and isothermal waterflooding strategies. The number of manipulated inputs, i.e. independent decision variables, is the same for both strategies. However, the thermal model contains more differential equations and algebraic equations.
The optimization of the thermal strategy involves 15 iterations in the SQLP optimization algorithm, whereas the optimization of the isothermal strategy only requires 6. However, the computation time per iteration is slightly higher for the isothermal strategy because each iteration, on average, involves more simulations. The number of gradient evaluations per iteration is close to 1 for both strategies.

The average number of time steps in the thermal and isothermal simulations are close to each other. However, the simulations of the thermal waterflooding strategies require, on average, close to one more Newton iteration. Consequently, the number of evaluations of the functions $F$ and $G$ in the semi-explicit differential-algebraic equations \eqref{eq:ocp:alg}-\eqref{eq:ocp:dyn}, and their Jacobians, is also approximately one higher. Furthermore, the GMRES algorithm requires close to 0.8 more iterations, on average, to solve the linear systems in the simulations of the thermal strategy. In conclusion, it is more computationally demanding to simulate, and therefore also optimize, the thermal waterflooding strategy than the isothermal strategy.
\begin{table}
	\centering
	\caption{Problem characteristics and KPIs for the dynamic optimization of the thermal and isothermal waterflooding strategies. The function evaluations refer to the evaluation of the functions $F$ and $G$ in the semi-explicit differential-algebraic equations \eqref{eq:ocp:alg}-\eqref{eq:ocp:dyn}, and the Jacobian evaluations refer to the evaluation of the Jacobians of $F$ and $G$. The iterations per linear system refer to the iterations in the GMRES algorithm. The KPIs related to the simulations are averaged over all the simulations involved in the dynamic optimization.}
	\label{tab:kpi}
	\begin{tabular}{lrr}
		& Thermal & Isothermal \\
		\cline{2-3}
		\textbf{Problem}                &        &        \\
		\hline
		Manipulated inputs              &    180 &    180 \\
		Differential equations          &    847 &    726 \\
		Algebraic equations             &   2541 &   2299 \\
		                                &        &        \\
		\textbf{Optimization}           &        &        \\
		\hline
		Iterations                      &     15 &      6 \\
		Simulations                     &     27 &     20 \\
		Gradient evaluations            &     16 &      7 \\
		CPU time (s)                    & 661.16 & 287.18 \\
		CPU time per iteration (s)      &  44.08 &  47.86 \\
		                                &        &        \\
		\textbf{Simulation}                     &        &        \\
		\hline
		Time steps per simulation          & 218.37 & 212.55 \\
		Newton iterations per time step    &   4.04 &   3.02 \\
		Function evaluations per time step &   6.07 &   5.07 \\
		Jacobian evaluations per time step &   5.67 &   4.41 \\
		Iterations per linear system       &   9.03 &   8.21 \\
	\end{tabular}
\end{table}

%% file: tex/concl.tex
\section{Conclusions}\label{sec:concl}
In this work, we consider dynamic optimization of thermal and isothermal oil recovery processes. Therefore, we present thermodynamically rigorous models of thermal and isothermal multicomponent three-phase flow in subsurface oil reservoirs.
The involved phase equilibrium problems, i.e. the UV and the VT flash, are based on the second law of thermodynamics. Furthermore, we formulate the UV and VT flash problems as optimization problems. We demonstrate that the thermal and the isothermal reservoir flow models are in a semi-explicit index-1 differential-algebraic form, and we use a gradient-based algorithm to solve the dynamic optimization problems. We implement the algorithm in C/C++ using the software libraries DUNE, ThermoLib, and KNITRO. Finally, we present numerical examples of optimized thermal and isothermal oil recovery strategies, and we discuss the computational performance of the dynamic optimization algorithm.

%% file: tex/app_relperm.tex
\section{Relative permeability}\label{app:relperm}
In this appendix, we use Stone's model II \cite{Delshad:Pope:1989} to describe the relative permeabilities. The model equations involve the relative permeabilities of 1) a hypothetical oil-water system and 2) a hypothetical oil-gas system. We express the relative permeabilities for these two hypothetical systems using the modified Brooks-Corey model \cite{Goda:Behrenbruch:2004}. We introduce the normalized saturations, $\bar{S}^w = \bar{S}^w(T, P, n^w, n^o, n^g)$ and $\bar{S}^g = \bar{S}^g(T, P, n^w, n^o, n^g)$:
\begin{subequations}
	\begin{align}
		\bar{S}^w &= \big(\hat{S}^w - \hat{S}_c^w\big)\big/\big(1 - \hat{S}_c^w - \hat{S}_{\max}^w\big), \\
		\bar{S}^g &= \big(\hat{S}^g - \hat{S}_c^g\big)\big/\big(1 - \hat{S}_c^g - \hat{S}_{\max}^g\big).
	\end{align}
\end{subequations}
$\hat{S}^w = \hat{S}^w(T, P, n^w, n^o, n^g)$ and $\hat{S}^g = \hat{S}^g(T, P, n^w, n^o, n^g)$ are the water and gas saturations, $\hat{S}_c^w$ and $\hat{S}_c^g$ are the connate water and gas saturations, and $\hat{S}_{\max}^w$ and $\hat{S}_{\max}^g$ are the maximum water and gas saturations. The relative permeabilities of the hypothetical oil-water system, $k_r^w = k_r^w(T, P, n^w, n^o, n^g)$ and $k_r^{ow} = k_r^{ow}(T, P, n^w, n^o, n^g)$, are
\begin{subequations}\label{eq:relperm}
	\begin{align}
		\label{eq:krw}
		k_r^w &= k_{r, 0}^w(\bar{S}^w)^{m_w}, \\
		\label{eq:krow}
		k_r^{ow} &= k_{r, 0}^{ow} \big(1 - \bar{S}^w\big)^{m_{ow}},
	\end{align}
	and similarly, the relative permeabilities of the hypothetical oil-gas system, $k_r^g = k_r^g(T, P, n^w, n^o, n^g)$ and $k_r^{og} = k_r^{og}(T, P, n^w, n^o, n^g)$, are
	\begin{align}
		\label{eq:krg}
		k_r^g &= k_{r, 0}^g(\bar{S}^g)^{m_g}, \\
		\label{eq:krog}
		k_r^{og} &= k_{r, 0}^{og} \big(1 - \bar{S}^g\big)^{m_{og}},
	\end{align}
\end{subequations}
where $k_{r, 0}^w$, $k_{r, 0}^{ow}$, $k_{r, 0}^g$, and $k_{r, 0}^{og}$ are the end-point relative permeabilities, and $m_w$, $m_{ow}$, $m_g$, and $m_{og}$ are the Corey exponents.
If one of the expressions in \eqref{eq:relperm} become negative or larger than one, the corresponding relative permeability is set to zero or one, respectively.
Finally, the relative permeability of the oil phase, $k_r^o = k_r^o(T, P, n^w, n^o, n^g)$, is
\begin{align}\label{eq:kro}
	k_r^o = k_r^c\big((k_r^{ow}/k_r^c + k_r^w)(k_r^{og}/k_r^c + k_r^g) - (k_r^w + k_r^g)\big),
\end{align}
where $k_r^c$ is a parameter. To summarize, the relative permeabilities of the water, oil, and gas phases are given by \eqref{eq:krw}, \eqref{eq:kro}, and \eqref{eq:krg}, respectively.

%% file: tex/app_visc.tex
\section{Viscosity of oil and gas}\label{app:visc}
In this appendix, we describe the model of the viscosity of reservoir fluids by Lohrenz et al. \cite{Lohrenz:etal:1964}. We use the expressions for the liquid phase viscosity to describe the viscosity of both the oil and the gas phase. The viscosity of phase $\alpha\in\{o, g\}$, $\mu^\alpha = \mu^\alpha(T, P, n^\alpha)$, is a function of temperature, $T$, pressure, $P$, and phase composition (in moles), $n^\alpha$:
\begin{align}
	\mu^\alpha &= \bar{\mu}^\alpha + \frac{1}{\tau^\alpha}\left(\left(a^\alpha\right)^4 - 10^{-4}\right), & \alpha &\in \{o, g\}.
\end{align}
The auxiliary variables $\tau^\alpha = \tau^\alpha(n^\alpha)$ and $a^\alpha = a^\alpha(T, P, n^\alpha)$ are
\begin{subequations}\label{eq:aux}
	\begin{align}
		\label{eq:aux:tau}
		\tau^\alpha &= \left(T_c^\alpha\right)^\frac{1}{6} \left(M_w^\alpha\right)^{-\frac{1}{2}} \left(P_c^\alpha\right)^{-\frac{2}{3}}, & \alpha &\in \{o, g\}, \\
		\label{eq:aux:a}
		a^\alpha &= \sum_{i=0}^4 a_i \left(\rho_r^\alpha\right)^i, & \alpha &\in \{o, g\}.
	\end{align}
\end{subequations}
The coefficients in the polynomial in \eqref{eq:aux:a} are $a_0 = 0.1023$, $a_1 = 0.023364$, $a_2 = 0.058533$, $a_3 = -0.040758$, and $a_4 = 0.0093324$ \cite{Jossi:etal:1962}. The auxiliary variables $T_c^\alpha = T_c^\alpha(n^\alpha)$, $P_c^\alpha = P_c^\alpha(n^\alpha)$, $V_c^\alpha = V_c^\alpha(n^\alpha)$, and $M_w^\alpha = M_w^\alpha(n^\alpha)$ are
\begin{subequations}\label{eq:phasecrit}
	\begin{align}
		T_c^\alpha &= \sum_{k=1}^{N_C} x_k^\alpha T_{c, k}, & \alpha &\in \{o, g\}, \\
		P_c^\alpha &= \sum_{k=1}^{N_C} x_k^\alpha P_{c, k}, & \alpha &\in \{o, g\}, \\
		V_c^\alpha &= \sum_{k=1}^{N_C} x_k^\alpha V_{c, k}, & \alpha &\in \{o, g\}, \\
		M_w^\alpha &= \sum_{k=1}^{N_C} x_k^\alpha M_{w, k}, & \alpha &\in \{o, g\}.
	\end{align}
\end{subequations}
We present the expression for $V_c^\alpha$ here, but we first use it in \eqref{eq:visc:dens:red}. We use values of the pure component critical temperature, $T_{c, k}$, critical pressure, $P_{c, k}$, critical volume, $V_{c, k}$, and molecular weight, $M_{w, k}$, from the DIPPR database \cite{Thomson:1996}. In order to describe the viscosity of the oil and the gas phases in a unified manner, we have adopted a different notation for the mole fractions, $x_k^\alpha = x_k^\alpha(n^\alpha)$, than in previous sections:
\begin{align}
	x_k^\alpha &= \frac{n_k^\alpha}{N^\alpha}, & k &= 1, \ldots, N_C, & \alpha &\in \{o, g\}. 
\end{align}
The total amount of moles in phase $\alpha$, $N^\alpha = N^\alpha(n^\alpha)$, is
\begin{align}
	N^\alpha &= \sum_{k=1}^{N_C} n_k^\alpha, & \alpha &\in \{o, g\}.
\end{align}
The reduced density, $\rho_r^\alpha = \rho_r^\alpha(T, P, n^\alpha)$, is
\begin{align}\label{eq:visc:dens:red}
	\rho_r^\alpha &= \rho^\alpha V_c^\alpha, & \alpha &\in \{o, g\},
\end{align}
where the molar density, $\rho^\alpha = \rho^\alpha(T, P, n^\alpha)$, is
\begin{align}
	\rho^\alpha &= \frac{N^\alpha}{V^\alpha}, & \alpha &\in \{o, g\}.
\end{align}
The reference viscosity, $\bar{\mu}^\alpha = \bar{\mu}^\alpha(T, n^\alpha)$, is
\begin{align}
	\bar{\mu}^\alpha &= \frac{\bar{M}_{w\mu}^\alpha}{\bar{M}_w^\alpha}, & \alpha &\in \{o, g\},
\end{align}
where the auxiliary variables $\bar{M}_{w\mu}^\alpha = \bar{M}_{w\mu}^\alpha(T, n^\alpha)$ and $\bar{M}_w^\alpha = \bar{M}_w^\alpha(n^\alpha)$ are
\begin{subequations}
	\begin{align}
		\bar{M}_{w\mu}^\alpha 	&= \sum_{k=1}^{N_C} x_k^\alpha \bar{\mu}_k \sqrt{M_{w, k}}, & \alpha &\in \{o, g\}, \\
		\bar{M}_w^\alpha 		&= \sum_{k=1}^{N_C} x_k^\alpha             \sqrt{M_{w, k}}, & \alpha &\in \{o, g\}.
	\end{align}
\end{subequations}
The pure component reference viscosity, $\bar{\mu}_k = \bar{\mu}_k(T)$, is
\begin{align}
	\bar{\mu}_k &=
	\left\{
	\begin{array}{ll}
		34\cdot 10^{-5} \frac{T_{r, k}^{0.94}}{\tau_k}, & T_{r, k} < 1.5, \\
		17.78\cdot 10^{-5}\frac{(4.58 T_{r, k} - 1.67)^\frac{5}{8}}{\tau_k}, & T_{r, k} \geq 1.5,
	\end{array}
	\right. & k &= 1, \ldots, N_C,
\end{align}
where the reduced temperature, $T_{r, k} = T_{r, k}(T)$, and the auxiliary variable $\tau_k$ are
\begin{subequations}
	\begin{align}
		T_{r, k} &= \frac{T}{T_{c, k}}, & k &= 1, \ldots, N_C, \\
		\tau_k &= T_{c, k}^\frac{1}{6} M_{w, k}^{-\frac{1}{2}} P_{c, k}^{-\frac{2}{3}}, & k &= 1, \ldots, N_C.
	\end{align}
\end{subequations}

%% file: ms.bbl
\begin{thebibliography}{10}

\bibitem{Alvarado:Manrique:2010}
{\sc V.~Alvarado and E.~Manrique}, {\em Enhanced oil recovery: an update
  review}, Energies, 3 (2010), pp.~1529--1575.

\bibitem{Bastian:etal:2008}
{\sc P.~Bastian, M.~Blatt, A.~Dedner, C.~Engwer, R.~Kl\"{o}fkorn, R.~Kornhuber,
  M.~Ohlberger, and O.~Sander}, {\em A generic grid interface for parallel and
  adaptive scientific computing. {P}art {II}: implementation and tests in
  {DUNE}}, Computing, 82 (2008), pp.~121--138.

\bibitem{Bastian:etal:2008b}
{\sc P.~Bastian, M.~Blatt, A.~Dedner, C.~Engwer, R.~Kl\"{o}fkorn, M.~Ohlberger,
  and O.~Sander}, {\em A generic grid interface for parallel and adaptive
  scientific computing. {P}art {I}: abstract framework}, Computing, 82 (2008),
  pp.~103--119.

\bibitem{Berning:etal:2009}
{\sc T.~Berning, M.~Odgaard, and S.~K. K{\ae}r}, {\em A computational analysis
  of multiphase flow through {PEMFC} cathode porous media using the multifluid
  approach}, Journal of the Electrochemical Society, 156 (2009),
  pp.~B1301--B1311.

\bibitem{Binder:etal:2001b}
{\sc T.~Binder, L.~Blank, H.~G. Bock, R.~Bulirsch, W.~Dahmen, M.~Diehl,
  T.~Kronseder, W.~Marquardt, J.~P. Schl{\"o}der, and O.~von Stryk}, {\em
  Introduction to model based optimization of chemical processes on moving
  horizons}, in Online Optimization of Large Scale Systems, Springer-Verlag
  Berlin Heidelberg, 2001, pp.~295--339.

\bibitem{Blatt:Bastian:2007}
{\sc M.~Blatt and P.~Bastian}, {\em The iterative solver template library}, in
  Applied Parallel Computing. State of the Art in Scientific Computing. PARA
  2006, {B. K{\aa}gstr\"{o}m et al.}, ed., vol.~4699 of Lecture Notes in
  Computer Science, Springer-Verlag Berlin Heidelberg, 2007, pp.~666--675.

\bibitem{Bukshtynov:etal:2015}
{\sc V.~Bukshtynov, O.~Volkov, L.~J. Durlofsky, and K.~Aziz}, {\em
  Comprehensive framework for gradient-based optimization in closed-loop
  reservoir management}, Computational Geosciences, 19 (2015), pp.~877--897.

\bibitem{Callen:1985}
{\sc H.~B. Callen}, {\em Thermodynamics and an introduction to
  thermostatistics}, John Wiley \& Sons, 2nd~ed., 1985.

\bibitem{Capolei:etal:2012}
{\sc A.~Capolei, C.~V{\"o}lcker, J.~Frydendall, and J.~B. J{\o}rgensen}, {\em
  Oil reservoir production optimization using single shooting and {ESDIRK}
  methods}, IFAC Proceedings Volumes, 45 (2012), pp.~286--291.

\bibitem{Chen:2007}
{\sc Z.~Chen}, {\em Reservoir simulation: mathematical techniques in oil
  recovery}, SIAM, 2007.

\bibitem{Chen:2006}
{\sc Z.~Chen, G.~Huan, and Y.~Ma}, {\em Computational methods for multiphase
  flows in porous media}, Computational Science \& Engineering, SIAM, 2006.

\bibitem{Codas:etal:2017}
{\sc A.~Codas, K.~G. Hanssen, B.~Foss, A.~Capolei, and J.~B. J{\o}rgensen},
  {\em Multiple shooting applied to robust reservoir control optimization
  including output constraints on coherent risk measures}, Computational
  Geosciences, 21 (2017), pp.~479--497.

\bibitem{Colson:etal:2007}
{\sc B.~Colson, P.~Marcotte, and G.~Savard}, {\em An overview of bilevel
  optimization}, Annals of Operations Research, 153 (2007), pp.~235--256.

\bibitem{Delshad:Pope:1989}
{\sc M.~Delshad and G.~A. Pope}, {\em Comparison of the three-phase oil
  relative permeability models}, Transport in Porous Media, 4 (1989),
  pp.~59--83.

\bibitem{Eppelbaum:etal:2014}
{\sc L.~Eppelbaum, I.~Kutasov, and A.~Pilchin}, {\em Applied geothermics},
  Lecture Notes in Earth System Sciences, Springer-Verlag Berlin Heidelberg,
  2014.

\bibitem{Forouzanfar:etal:2013}
{\sc F.~Forouzanfar, E.~D. Rossa, R.~Russo, and A.~C. Reynolds}, {\em
  Life-cycle production optimization of an oil field with an adjoint-based
  gradient approach}, Journal of Petroleum Science and Engineering, 112 (2013),
  pp.~351--358.

\bibitem{Garipov:etal:2018}
{\sc T.~T. Garipov, P.~Tomin, R.~Rin, D.~V. Voskov, and H.~A. Tchelepi}, {\em
  Unified thermo-compositional-mechanical framework for reservoir simulation},
  Computational Geosciences,  (2018),
  \url{https://doi.org/10.1007/s10596-018-9737-5}.

\bibitem{Gmehling:etal:2012}
{\sc J.~Gmehling, B.~Kolbe, M.~Kleiber, and J.~Rarey}, {\em Chemical
  thermodynamics for process simulation}, Wiley-VCH, 2012.

\bibitem{Goda:Behrenbruch:2004}
{\sc H.~M. Goda and P.~Behrenbruch}, {\em Using a modified {B}rooks-{C}orey
  model to study oil-water relative permeability for diverse pore structures},
  in Proceedings of the 2004 SPE Asia Pacific Oil and Gas Conference and
  Exhibition, Perth, Australia, Oct. 2004.

\bibitem{Hannaoui:etal:2015}
{\sc R.~Hannaoui, P.~Horgue, F.~Larachi, Y.~Haroun, F.~Augier, M.~Quintard, and
  M.~Prat}, {\em Pore-network modeling of trickle bed reactors: pressure drop
  analysis}, Chemical Engineering Journal, 262 (2015), pp.~334--343.

\bibitem{Heirung:etal:2011}
{\sc T.~A.~N. Heirung, M.~R. Wartmann, J.~D. Jansen, B.~E. Ydstie, and B.~A.
  Foss}, {\em Optimization of the water-flooding process in a small {2D}
  horizontal oil reservoir by direct transcription}, IFAC Proceedings Volumes,
  44 (2011), pp.~10863--10868.

\bibitem{Holman:2010}
{\sc J.~P. Holman}, {\em Heat transfer}, McGraw-Hill, 10th~ed., 2010.

\bibitem{Jang:Aral:2009}
{\sc W.~Jang and M.~M. Aral}, {\em Multiphase flow fields in in-situ air
  sparging and its effect on remediation}, Transport in Porous Media, 76
  (2009), pp.~99--119.

\bibitem{Jorgensen:2007b}
{\sc J.~B. J{\o}rgensen}, {\em Adjoint sensitivity results for predictive
  control, state- and parameter-estimation with nonlinear models}, in
  Proceedings of the 2007 European Control Conference, Kos, Greece, July 2007,
  pp.~3649--3656.

\bibitem{Jossi:etal:1962}
{\sc J.~A. Jossi, L.~I. Stiel, and G.~Thodos}, {\em The viscosity of pure
  substances in the dense gaseous and liquid phases}, A.I.Ch.E Journal, 8
  (1962), pp.~59--63.

\bibitem{Khan:etal:2018}
{\sc M.~I.~H. Khan, M.~U.~H. Joardder, C.~Kumar, and M.~A. Karim}, {\em
  Multiphase porous media modelling: a novel approach to predicting food
  processing performance}, Critical Reviews in Food Science and Nutrition, 58
  (2018), pp.~528--546.

\bibitem{Kone:etal:2017}
{\sc J.-P. Kone, X.~Zhang, Y.~Yan, G.~Hu, and G.~Ahmadi}, {\em
  Three-dimensional multiphase flow computational fluid dynamics models for
  proton exchange membrane fuel cell: a theoretical development}, The Journal
  of Computational Multiphase Flows, 9 (2017), pp.~3--25.

\bibitem{Koretsky:2013}
{\sc M.~D. Koretsky}, {\em Engineering and chemical thermodynamics}, Wiley,
  2nd~ed., 2013.

\bibitem{Kourounis:etal:2014}
{\sc D.~Kourounis, L.~J. Durlofsky, J.~D. Jansen, and K.~Aziz}, {\em Adjoint
  formulation and constraint handling for gradient-based optimization of
  compositional reservoir flow}, Computational Geosciences, 18 (2014),
  pp.~117--137.

\bibitem{Kristensen:etal:2005}
{\sc M.~R. Kristensen, J.~B. J{\o}rgensen, P.~G. Thomsen, M.~L. Michelsen, and
  S.~B. J{\o}rgensen}, {\em Sensitivity analysis in index-1 differential
  algebraic equations by {ESDIRK} methods}, IFAC Proceedings Volumes, 38
  (2005), pp.~212--217.

\bibitem{Lei:etal:2012}
{\sc Y.~Lei, S.~Li, X.~Zhang, and Q.~Zhang}, {\em Optimal control of polymer
  flooding for high temperature and high salinity reservoir}, International
  Journal of Advancements in Computing Technology, 4 (2012), pp.~52--60.

\bibitem{Lie:2014}
{\sc K.-A. Lie}, {\em An introduction to reservoir simulation using {MATLAB}},
  Sintef ICT, Oslo, Norway, 2014.

\bibitem{Lohrenz:etal:1964}
{\sc J.~Lohrenz, B.~G. Bray, and C.~R. Clark}, {\em Calculating viscosities of
  reservoir fluids from their compositions}, Journal of Petroleum Technology,
  16 (1964), pp.~1171--1176.

\bibitem{Luo:etal:1996}
{\sc Z.-Q. Luo, J.-S. Pang, and D.~Ralph}, {\em Mathematical programs with
  equilibrium constraints}, Cambridge University Press, 1996.

\bibitem{Michelsen:1999}
{\sc M.~L. Michelsen}, {\em State function based flash specifications}, Fluid
  Phase Equilibria, 158-160 (1999), pp.~617--626.

\bibitem{Nocedal:Wright:2006}
{\sc J.~Nocedal and S.~J. Wright}, {\em Numerical optimization}, Springer
  Science \& Business Media, 2nd~ed., 2006.

\bibitem{Onwunalu:Durlofsky:2010}
{\sc J.~E. Onwunalu and L.~J. Durlofsky}, {\em Application of a particle swarm
  optimization algorithm for determining optimum well location and type},
  Computational Geosciences, 14 (2010), pp.~183--198.

\bibitem{Outrata:etal:1998}
{\sc J.~Outrata, M.~Ko{\v{c}}vara, and J.~Zowe}, {\em Nonsmooth approach to
  optimization problems with equilibrium constraints: theory, applications and
  numerical results}, vol.~28 of Nonconvex Optimization and Its Applications,
  Springer Science \& Business Media, 1998.

\bibitem{Peng:robinson:1976}
{\sc D.-Y. Peng and D.~B. Robinson}, {\em A new two-constant equation of
  state}, Industrial \& Engineering Chemistry Fundamentals, 15 (1976),
  pp.~59--64.

\bibitem{Pesavento:etal:2017}
{\sc F.~Pesavento, B.~A. Schrefler, and G.~Scium{\`{e}}}, {\em Multiphase flow
  in deforming porous media: a review}, Archives of Computational Methods in
  Engineering, 24 (2017), pp.~423--448.

\bibitem{Polivka:Mikyska:2014}
{\sc O.~Pol{\'{i}}vka and J.~Miky{\v{s}}ka}, {\em Compositional modeling in
  porous media using constant volume flash and flux computation without the
  need for phase identification}, Journal of Computational Physics, 272 (2014),
  pp.~149--169.

\bibitem{Ritschel:etal:2017}
{\sc T.~K.~S. Ritschel, A.~Capolei, J.~Gaspar, and J.~B. J{\o}rgensen}, {\em An
  algorithm for gradient-based dynamic optimization of {UV} flash processes},
  Computers and Chemical Engineering, 114 (2018), pp.~281--295.

\bibitem{Ritschel:Capolei:Jorgensen:2017b}
{\sc T.~K.~S. Ritschel, A.~Capolei, and J.~B. J{\o}rgensen}, {\em The adjoint
  method for gradient-based dynamic optimization of {UV} flash processes},
  Computer Aided Chemical Engineering, 40 (2017), pp.~2071--2076.

\bibitem{Ritschel:Capolei:Jorgensen:2017}
{\sc T.~K.~S. Ritschel, A.~Capolei, and J.~B. J{\o}rgensen}, {\em Dynamic
  optimization of {UV} flash processes}, in FOCAPO / CPC 2017, Tucson, Arizona,
  Jan. 2017.

\bibitem{Ritschel:etal:2016}
{\sc T.~K.~S. Ritschel, J.~Gaspar, A.~Capolei, and J.~B. J{\o}rgensen}, {\em An
  open-source thermodynamic software library}, Tech. Report DTU Compute
  Technical Report-2016-12, Department of Applied Mathematics and Computer
  Science, Technical University of Denmark, 2016.

\bibitem{Ritschel:etal:2017b}
{\sc T.~K.~S. Ritschel, J.~Gaspar, and J.~B. J{\o}rgensen}, {\em A
  thermodynamic library for simulation and optimization of dynamic processes},
  IFAC PapersOnLine, 50 (2017), pp.~3542--3547.

\bibitem{Ritschel:Jorgensen:2017}
{\sc T.~K.~S. Ritschel and J.~B. J{\o}rgensen}, {\em Computation of phase
  equilibrium and phase envelopes}, Tech. Report DTU Compute Technical
  Report-2017-11, Department of Applied Mathematics and Computer Science,
  Technical University of Denmark, 2017.

\bibitem{Ritschel:Jorgensen:2018b}
{\sc T.~K.~S. Ritschel and J.~B. J{\o}rgensen}, {\em Computation of phase
  equilibrium in reservoir simulation and optimization}, in Proceedings of the
  3rd IFAC Workshop on Automatic Control in Offshore Oil and Gas Production,
  Esbjerg, Denmark, May 2018.

\bibitem{Ritschel:Jorgensen:2018c}
{\sc T.~K.~S. Ritschel and J.~B. J{\o}rgensen}, {\em Production optimization of
  a rigorous thermal and compositional reservoir flow model}, in Proceedings of
  the 3rd IFAC Workshop on Automatic Control in Offshore Oil and Gas
  Production, Esbjerg, Denmark, May 2018.

\bibitem{Ritschel:Jorgensen:2018e}
{\sc T.~K.~S. Ritschel and J.~B. J{\o}rgensen}, {\em Production optimization of
  thermodynamically rigorous isothermal and compositional models}, in
  Proceedings of the 16th European Conference on the Mathematics of Oil
  Recovery, Barcelona, Spain, Sept. 2018.

\bibitem{Simon:Ulbrich:2015}
{\sc M.~Simon and M.~Ulbrich}, {\em Adjoint based optimal control of partially
  miscible two-phase flow in porous media with applications to {CO$_2$}
  sequestration in underground reservoirs}, Optimization and Engineering, 16
  (2015), pp.~103--130.

\bibitem{Smith:etal:2005}
{\sc J.~M. Smith, H.~C. {Van Ness}, and M.~M. Abbott}, {\em Introduction to
  chemical engineering thermodynamics}, McGraw-Hill, 7th~ed., 2005.

\bibitem{Thomson:1996}
{\sc G.~H. Thomson}, {\em The {DIPPR}{\textregistered} databases},
  International Journal of Thermophysics, 17 (1996), pp.~223--232.

\bibitem{Thorvaldsson:Palsson:2012}
{\sc L.~Thorvaldsson and H.~Palsson}, {\em Modeling liquid dominated two phase
  flow in geothermal reservoirs in vicinity to, and inside wells}, in
  Proceedings of the 37th Workshop on Geothermal Reservoir Engineering,
  Stanford, California, USA, Jan. 2012.

\bibitem{Volcker:etal:2011}
{\sc C.~V{\"o}lcker, J.~B. J{\o}rgensen, and E.~H. Stenby}, {\em Oil reservoir
  production optimization using optimal control}, in Proceedings of the 50th
  {IEEE} Conference on Decision and Control and European Control Conference,
  Orlando, Florida, USA, Dec. 2011, pp.~7937--7943.

\bibitem{Volcker:etal:2010}
{\sc C.~V{\"o}lcker, J.~B. J{\o}rgensen, P.~G. Thomsen, and E.~H. Stenby}, {\em
  Adaptive stepsize control in implicit {R}unge-{K}utta methods for reservoir
  simulation}, IFAC Proceedings Volumes, 43 (2010), pp.~523--528.

\bibitem{Wachter:Biegler:2006}
{\sc A.~W{\"a}chter and L.~T. Biegler}, {\em On the implementation of an
  interior-point filter line-search algorithm for large-scale nonlinear
  programming}, Mathematical Programming, 106 (2006), pp.~25--57.

\bibitem{Zhang:Li:2007}
{\sc Z.~Xiao-Dong and L.~Shu-rong}, {\em Optimal control solving of polymer
  flooding in enhanced oil recovery with {2-D} models}, in Proceedings of the
  2007 IEEE International Conference on Control and Automation, Guangzhou,
  China, May 2007, pp.~1981--1986.

\bibitem{Zaydullin:etal:2014}
{\sc R.~Zaydullin, D.~V. Voskov, S.~C. James, H.~Henley, and A.~Lucia}, {\em
  Fully compositional and thermal reservoir simulation}, Computers and Chemical
  Engineering, 63 (2014), pp.~51--65.

\bibitem{Zhao:etal:2016}
{\sc H.~Zhao, Y.~W. Tang, Y.~Li, Y.~B. Shi, L.~Cao, R.~X. Gong, and G.~H.
  Shang}, {\em Reservoir production optimization using general stochastic
  approximate algorithm under the mixed-linear-nonlinear constraints}, Journal
  of Residuals Science \& Technology, 13 (2016).

\end{thebibliography}
